

\input amstex
\documentstyle{amsppt}

\input label.def
\input degt.def
\def\mnote#1{}

\input epsf
\def\picture#1{\epsffile{#1-bb.eps}}

\def\ie.{\emph{i.e\.}}
\def\eg.{\emph{e.g\.}}
\def\cf.{\emph{cf\.}}
\def\via{\emph{via}}

\def\Apriori{\emph{A priori}}

{\catcode`\@11
\gdef\proclaimfont@{\sl}
\gdef\subsubheadfont@{\subheadfont@}
}

\newtheorem{Remark}\Remark\thm\endAmSdef
\problem\thm\endAmSdef
\conjecture\thm\endproclaim
\newhead\subsubsection\subsubsection\endsubhead
\def\paragraph{\subsubsection{}}

\def\dash{\item"\hfill--\hfill"}
\def\Dashes{\widestnumber\item{--}\roster}
\def\endDashes{\endroster}

\loadbold
\def\bA{\bold A}
\def\bD{\bold D}
\def\bE{\bold E}

\def\tA{\tilde\bA}
\def\tD{\tilde\bD}
\def\tE{\tilde\bE}

\let\Ga\alpha
\let\Gb\beta
\let\Gg\gamma
\let\Gd\delta
\let\Gs\sigma
\let\Ge\varepsilon
\let\Gf\varphi
\let\Gr\rho

\def\dR{\partial R}

\def\Cp#1{\Bbb P^{#1}}

\def\Pic{\operatorname{Pic}}
\def\discr{\operatorname{discr}}
\def\Rad{\operatorname{ker}}
\def\res{\operatorname{res}}

\let\tree\Xi

\let\<\langle
\let\>\rangle
\def\1{^{-1}}
\def\BG#1{\Bbb B_{#1}}
\def\CG#1{\Z_{#1}}

\def\PSL{\operatorname{\text{\sl PSL}}}
\def\SL{\operatorname{\text{\sl SL}}}
\def\O{\operatorname{\text{\sl O}}}
\def\MW{\operatorname{\text{\sl MW}}}
\def\NS{\operatorname{\text{\sl NS}}}
\def\op{\operatorname{op}}
\def\nx{\operatorname{nx}}
\def\m#1{\bar#1}

\def\CE{\Cal E}
\let\lattice\Cal
\def\CH{\lattice H}
\def\CL{\lattice L}
\def\CU{\lattice U}
\def\CV{\lattice V}
\def\CW{\lattice W}
\def\CQ{\lattice Q}
\def\CT{\lattice T}
\def\CS{\lattice S}
\def\be{\bold e}

\def\bq{\bold q}

\def\bu{\bold u}
\def\bv{\bold v}
\def\bw{\bold w}
\def\bx{\bold x}

\def\Xcirc{X^\circ}
\def\Fcirc{F^\circ}
\def\Bcirc{B^\circ}
\def\dX{\partial\Xcirc}
\def\barF{\bar F}
\def\bz{\bar z}
\def\X{\Bbb X}
\def\Y{\Bbb Y}

\def\mm{\frak m}
\def\mh{\frak h}
\def\pp#1{\mathopen\|#1\mathclose\|}
\def\fc#1{[#1]}
\let\into\hookrightarrow

\def\DD{
\def\bul##1{\hbox to0pt{\hss$\overset##1\to\bullet$\hss}}
\def\-{\joinrel\relbar\joinrel\relbar\joinrel\relbar\joinrel\relbar\joinrel}
\def\={\joinrel\relbar\joinrel\relbar\joinrel}
}

\let\contraction\rightarrowtail
\let\contraction\rightsquigarrow
\let\gtimes\otimes

\def\inserthyphen{\ifcat\next a-\fi\ignorespaces}
\let\BLACK\bullet
\let\WHITE\circ
\def\CROSS{\vcenter{\hbox{$\scriptstyle\mathord\times$}}}
\let\STAR*
\def\pblack-{$\BLACK$\futurelet\next\inserthyphen}
\def\pwhite-{$\WHITE$\futurelet\next\inserthyphen}
\def\pcross-{$\CROSS$\futurelet\next\inserthyphen}
\def\pstar-{$\STAR$\futurelet\next\inserthyphen}
\def\black{\protect\pblack}
\def\white{\protect\pwhite}

\def\edge{{\bullet}{\joinrel\relbar\joinrel\relbar\joinrel}{\circ}}
\def\oedge{{\bullet}{\joinrel\relbar\joinrel\leftarrow\joinrel}{\circ}}

\topmatter

\author
Alex Degtyarev
\endauthor

\title
Transcendental lattice of an extremal elliptic surface
\endtitle

\address
Department of Mathematics,
Bilkent University,
06800 Ankara, Turkey
\endaddress

\email
degt\@fen.bilkent.edu.tr
\endemail

\abstract
We develop an algorithm computing the transcendental lattice and
the Mordell--Weil group of an extremal elliptic surface. As an
example, we compute the lattices of four exponentially large
series of surfaces
\endabstract



\endtopmatter

\document

\section{Introduction\label{S.intro}}

\subsection{Principal results}\label{s.results}
An \emph{extremal elliptic surface} can be defined
as a Jacobian elliptic
surface~$X$ of maximal Picard number,
$\rank\NS(X)=h^{1,1}(X)$,\mnote{$\Pic\to\NS$}
and minimal Mordell-Weil rank, $\rank\MW(X)=0$. For alternative,
more topological descriptions, see
Definition~\ref{def.extremal} and
Remark~\ref{rem.positive}.

Extremal elliptic surfaces are rigid; they are defined over
algebraic number fields. Up to isomorphism, such a surface~$X$
(without type~$\tE$ singular fibers)
is determined by an oriented $3$-regular ribbon
graph~$\Gamma_X$, called \emph{skeleton} of~$X$, see
Subsection~\ref{s.GammaX}.
This intuitive approach gives one a simple way to construct and
classify extremal
elliptic\mnote{`elliptic'}
surfaces, see, \eg.,~\cite{Beukers}
or~\cite{degt.kplets}; however, the relation between the
invariants of~$X$ and the
structure of~$\Gamma_X$ is not yet well understood.
A few first attempts to compute the invariants of some surfaces
were recently made in~\cite{ArimaShimada}.
In\mnote{new text}
slightly different terms, general properties of the
(necessarily finite)
Mordell-Weil group of an extremal elliptic surface and a few
examples are found in~\cite{Shioda.modular}.
(Due to~\cite{Shioda} and Nikulin's theory of lattice
extensions~\cite{Nikulin}, the Mordell-Weil group and the
transcendental lattice are closely related,
\cf.~\ref{s.hom.type}.)

The principal results of this paper are
Theorem~\ref{tripods} and Corollaries~\ref{cor.CT}
and~\ref{cor.MW}, computing the transcendental lattice~$\CT_X$ and
the Mordell--Weil group $\MW(X)$ of an extremal elliptic surface~$X$
without type~$\tE$ singular fibers in terms of its
skeleton~$\Gamma_X$. (Some
generalizations to wider classes of surfaces
are discussed in Section~\ref{S.generalizations},
see Theorems~\ref{tripods.E} and~\ref{tripods.ne}.)
It is important to notice that the algorithm uses a computer
friendly presentation of the graph (by a pair of permutations,
see Remark~\ref{rem.computer}); combined with the known
classification results (see, \eg.,~\cite{Beukers}) and various
lattice analyzing software, it can be used for computer
experiments.

\subsection{Examples}\label{s.examples}
Originally, this paper was motivated by a construction
in~\cite{degt.kplets}, producing exponentially large series of
non-isomorphic
extremal elliptic surfaces. Here, we compute the invariants of
these surfaces.

Given an integer $k\ge1$, define the lattices
(see Subsection~\ref{s.lattices}) $\CV_{k-1}$ and
$\CW_k$ as the orthogonal
direct\mnote{`direct'}
sums
$$
\CV_{k-1}=\bigoplus_{i=1}^{k-1}\Z\bv_i,\qquad
\CW_k=\bigoplus_{i=1}^{k-1}\Z\bv_i\oplus\Z\bw,
\eqtag\label{eq.UV}
$$
where $\bv_i^2=1$, $i=1,\ldots,k-1$, and $\bw^2=0$.


\theorem\label{th.tree.0.even}
Let~$X$ be
an extremal\mnote{`extremal'}
elliptic surface with
singular fibers
$$
\tA_{10s-2}\oplus(2s+1)\tA_0^*,\quad s\ge1.
$$
Then
$\CT_X\cong(3\bv_1+\ldots+3\bv_s+\bv_{s+1}+\ldots+\bv_{2s-1})^\perp\subset\CV_{2s-1}$.
\endtheorem

\theorem\label{th.tree.0.odd}
Let~$X$ be
an extremal\mnote{`extremal'}
elliptic surface with
singular fibers
$$
\tD_{10s-2}\oplus(2s)\tA_0^*,\quad s\ge1.
$$
Then $\CT_X\cong\bD_{2s-2}$ \rom(where
we let
$\bD_0=0$, $\bD_1=[4]$, $\bD_2=2\bA_1$, and
$\bD_3=\bA_3$\rom).
\endtheorem

Theorems~\ref{th.tree.0.even} and~\ref{th.tree.0.odd}
are proved in Subsection~\ref{pf.tree.0}.

\theorem\label{th.tree.1.even}
Let~$X$ be
an extremal\mnote{`extremal'}
elliptic surface with
singular fibers
$$
\tD_{10s+3}\oplus\tD_5\oplus(2s)\tA_0^*,\quad s\ge1.
$$
Then $\CT_X\cong\bD_{2s-1}\oplus\Z\bx$, where $\bx^2=4$.
\endtheorem

Let
$f_s=3\bv_1+\ldots+3\bv_{s-1}+\bv_s+\ldots+\bv_{2s-2}\in\CV_{2s-2}$,
and denote by $\CV'_{2s-2}$ the group~$\CV_{2s-2}$ with the
bilinear form
$x\otimes y\mapsto x\cdot y+\frac14(f_s\cdot x)(f_s\cdot y)$,
where $\,\cdot\,$ stands for the original product in~$\CV_{2s-2}$.
(Certainly, $\CV'_{2s-2}$ is \emph{not} an integral lattice.)

\theorem\label{th.tree.1.odd}
Let~$X$ be
an extremal\mnote{`extremal'}
elliptic surface with
singular fibers
$$
\tA_{10s-7}\oplus\tD_5\oplus(2s-1)\tA_0^*,\quad s\ge1.
$$
Then
$\CT_X$ is the index~$4$ sublattice
$\{x\in\CV'_{2s-2}\,|\,f_s\cdot x=0\bmod4\}\subset\CV'_{2s-2}$.
\endtheorem

Theorems~\ref{th.tree.1.even} and~\ref{th.tree.1.odd}
are proved in Subsection~\ref{pf.tree.1}.

Note that, in Theorems~\ref{th.tree.0.even}--\ref{th.tree.1.odd},
a simple count using the Riemann--Hurwitz formula for the
$j$-invariant
shows that the base of any extremal elliptic surface with one of
the combinatorial types of singular fibers indicated in the
statements
is~$\Cp1$.\mnote{a remark added}

The Jacobian elliptic surfaces as in
Theorems~\ref{th.tree.0.even}--\ref{th.tree.1.odd} appeared
in~\cite{degt.kplets}; within each of the four series, the number
of fiberwise equisingular deformation classes grows faster than
$a^{4s}$ for any $a<2$, \cf.~\ref{s.marking}, and the original
goal of this project was to distinguish these surfaces
topologically, hoping that the definite lattices~$\CT_X$ would
fall into distinct isomorphism classes. The four theorems above
show that this approach fails.
(Note that the theorems imply as well that, for each surface~$X$
in question, the Mordell--Weil group $\MW(X)$ is trivial.)
To add to the disappointment, one
can also use~\cite{degt.kplets} and
some intermediate results of this paper and compute the fundamental
groups $\pi_1(\Sigma\sminus(C\cup E))$ of the ramification loci of
the double coverings $X\to\Sigma$, see~\ref{s.Sigma}. Most groups
turn out to be abelian; hence they
also depend on~$s$ only (within each
of the four series).

\theorem\label{th.group}
Let $X$ be one of the surfaces as in
Theorems~\ref{th.tree.0.even}--\ref{th.tree.1.odd}, and assume that
$s>1$. Then  the fundamental group
$\pi_1(\Sigma\sminus(C\cup E))$ is cyclic.
\endtheorem

This theorem is proved in Subsection~\ref{pf.group}. In the four
exceptional cases corresponding to
the value $s=1$, the groups can also be
computed; they are listed in Remark~\ref{rem.group}. In two cases,
the trigonal curve~$C$ is reducible.

Thus, neither~$\CT_X$ nor $\pi_1(\Sigma\sminus(C\cup E))$
distinguish the surfaces, and the following problem, which
motivated this paper, still stands.

\problem
Are surfaces~$X$ as in
Theorems~\ref{th.tree.0.even}--\ref{th.tree.1.odd}
fiberwise homeomorphic (for each given~$s$ and within each given
series)? Are they Galois conjugate?
\endproblem

An answer to the first question should be given by the Hurwitz
equivalence class of the braid monodromy of the ramification
locus. The monodromies are given by~\eqref{eq.braids};
at present, I
do not know whether they are Hurwitz equivalent.

\subsection{Contents of the paper}\label{s.contents}
In Section~\ref{S.prelim} we remind a few concepts related to
integral lattices and elliptic surfaces. Section~\ref{S.topology}
deals with the topological part of the computation; it is used in
the proof of the main theorem and its corollaries
in Section~\ref{S.main}. In Section~\ref{S.trees}, we consider a
special class of skeletons, the so called
\emph{pseudo-\allowbreak trees},
and prove Theorems~\ref{th.tree.0.even}--\ref{th.group}. Finally,
in Section~\ref{S.generalizations}, we discuss a few
generalizations of the principal results.

\subsection{Acknowledgements}
I am grateful to I.~Itenberg, with whom I
discussed this project at its early stages, and to
I.~Dolgachev\mnote{more acknowledgements} and I.~Shimada
for their helpful remarks. I would also like to extend my
gratitude to the referee of this paper, whose remarks helped me to
improve clarity of the exposition.

\section{Preliminaries\label{S.prelim}}

\subsection{Lattices}\label{s.lattices}
An \emph{\rom(integral\rom) lattice} is a
finitely generated
free abelian
group~$\CL$
supplied with a symmetric bilinear form
$\CL\otimes\CL\to\Z$ (which is usually referred to as \emph{product}
and denoted by $x\otimes y\mapsto x\cdot y$ and
$x\otimes x\mapsto x^2$).
A lattice is called \emph{even} if $x^2=0\bmod2$ for all
$x\in\CL$.
Occasionally, we will also consider
\emph{rational lattices}, which are free abelian groups supplied
with $\Q$-valued symmetric bilinear forms.
A lattice structure on~$\CL$
is uniquely determined by the function
$x\mapsto x^2$: one has $x\cdot y=\frac12[(x+y)^2-x^2-y^2]$.

Given a lattice~$\CL$, one can define the \emph{associated
homomorphism} $\Gf_{\CL}\:\CL\to\CL^*:=\Hom(\CL,\Z)$ \via\
$x\mapsto[y\mapsto x\cdot y]\in\CL^*$.
The \emph{kernel} $\Rad\CL$ is the
kernel of~$\Gf_{\CL}$.
(We use the notation $\Rad\CL$ for the kernel of a lattice
as opposed to $\Ker\Ga$ for the kernel of a homomorphism~$\Ga$.)
A lattice~$\CL$
is called \emph{nondegenerate} if $\Rad\CL=0$; it is called
\emph{unimodular} if $\Gf_{\CL}$ is an isomorphism. For example,
the \emph{intersection lattice} $H_2(X)/\Tors$ of an oriented
closed $4$-manifold~$X$ is unimodular (Poincar\'e duality).

We will fix the notation~$\CU$ for the \emph{hyperbolic plane},
which is the unimodular lattice generated by two elements $\bu_1$,
$\bu_2$ with $\bu_1^2=\bu_2^2=0$, $\bu_1\cdot\bu_2=1$. We will
also use the notation $\bA_p$, $\bD_q$, $\bE_6$, $\bE_7$, $\bE_8$
for the irreducible positive definite lattices generated by the
root systems of the same name.

\paragraph\label{s.discr}
If $\CL$ is nondegenerate, the
quotient $\discr\CL:=\CL^*\!/\CL$ is a finite group; it is called
the \emph{discriminant group} of~$\CL$.
Since $\Gf_{\CL}\otimes\Q$ is an
isomorphism, $\CL^*$ turns into a rational lattice and $\discr\CL$
inherits a $(\Q/\Z)$-valued symmetric bilinear form
$$
(x\bmod\CL)\otimes(y\bmod\CL)\mapsto(x\cdot y)\bmod\Z,
$$
called the \emph{discriminant form} of~$\CL$. In general,
if $\CL$ is degenerate, we
define $\discr\CL$ do be $\discr(\CL/\Rad\CL)$. As a group,
$\discr\CL=\Tors(\CL^*\!/\CL)$.

If $\CL$ is even, $\discr\CL$ inherits also a $(\Q/2\Z)$-valued
quadratic extension of the discriminant form; it is given by
$(x\bmod\CL)\mapsto x^2\bmod2\Z$.

\paragraph\label{s.extension}
Let~$\CL$ be a unimodular lattice, and let $\CS\subset\CL$ be a
nondegenerate primitive sublattice. Denote $\CT=\CS^\perp$; it is
also nondegenerate.
According to Nikulin~\cite{Nikulin}, the image of the restriction
homomorphism $\CL^*\to\CS^*\oplus\CT^*\to\discr\CS\oplus\discr\CT$
is the graph of a certain anti-isometry
$q\:\discr\CS\to\discr\CT$. (If $\CL$ is even, then so are $\CS$
and~$\CT$ and $q$ is also an anti-isometry of the quadratic
extensions.) Furthermore, the pair $(\CT,q)$, up to the action of
$\O(\CT)$ on $\discr\CT$, determines the isomorphism class of the
extension $\CL\supset\CS$.

\subsection{Elliptic surfaces}\label{s.surfaces}
Here, we remind a few facts concerning elliptic surfaces. The
references are~\cite{FM} or the original paper~\cite{Kodaira}.

A \emph{Jacobian elliptic surface} is a compact complex
surface $X$ equipped with an elliptic fibration $\pr\:X\to B$
(\ie., a fibration with all but finitely many fibers nonsingular
elliptic curves)
and a distinguished section $E\subset X$ of~$\pr$.
(From the existence of a section it follows that $X$ has no
multiple fibers.) Throughout the paper we assume that surfaces are
\emph{relatively minimal}, \ie., fibers of~$\pr$ contain no
$(-1)$-curves.

For the topological type of a singular elliptic fiber~$F$, we use
the notation~$\tA$, $\tD$, $\tE$ referring to the extended Dynkin
graph representing the adjacencies of the components of~$F$. The
advantage of this approach is the fact that it reflects the type
of the corresponding
singular point of the ramification locus of~$X$,
\cf.~\ref{s.Sigma}. For the relation to Kodaira's notation
$\fam0I$--$\fam0IV^*$, values of the $j$-invariant,
and some other invariants, see
Table~1 in~\cite{degt.kplets}.\mnote{ref added}

\paragraph\label{s.j}
Let $\Bcirc\subset B$ be the set of regular values of~$\pr$, and
define the \emph{\rom(functional\rom) $j$-invariant}
$j_X\:B\to\Cp1$ as the analytic
continuation of the function $\Bcirc\to\C^1$ sending each
nonsingular fiber to its classical $j$-invariant (divided
by~$12^3$).

The monodromy $\mh_X\:\pi_1(\Bcirc)\to\SL(2,\Z)$
(in the $1$-homology of the fiber)
of the locally trivial fibration
$\pr\1\Bcirc\to\Bcirc$ is called the \emph{homological invariant}
of~$X$. Its reduction to
$\PSL(2,\Z)=\SL(2,\Z)/\{\pm1\}$\mnote{`$/\{\pm1\}$'}
is determined by
the $j$-invariant.
Together, $j_X$ and~$\mh_X$ determine~$X$
up to isomorphism; conversely, any pair
$(j,\mh)$ that agrees in the sense just described gives rise to a
Jacobian elliptic surface.

In particular, the homological invariant determines the \emph{type
specification} of~$X$, \ie., a choice of type, $\tA$ or~$\tD$,
$\tE$, of each singular fiber. If the base~$B$ is rational, then
the type specification and~$j_X$ determine~$\mh_X$.

\definition\label{def.extremal}
A Jacobian elliptic surface~$X$ is called \emph{extremal} if it
satisfies the following conditions:
\roster
\item\local{ext.1}
$j_X$ has no critical values other than~$0$, $1$, and~$\infty$;
\item\local{ext.2}
each point in $j_X\1(0)$ has ramification index at most~$3$,
and each point in $j_X\1(1)$ has ramification index at most~$2$;
\item\local{ext.3}
$X$ has no singular fibers of types~$\tD_4$, $\tA_0^{**}$,
$\tA_1^*$, or~$\tA_2^*$.
\endroster
(In fact, this more topological definition is the contents
of~\cite{MNori}.)
\enddefinition

\paragraph\label{s.hom.type}
Let $X$ be a Jacobian elliptic surface. Denote by
$\Gs_X\subset H_2(X)$ the set of classes realized by the
components of the singular fibers of~$X$. (We assume that $X$ does
have at least one singular fiber.) Let $\CS_X\subset H_2(X)$ be
the sublattice spanned by~$\Gs_X$ and~$[E]$ (sometimes, $\CS_X$ is
called the \emph{simple lattice} of~$X$), and let
$\tilde\CS_X:=(\CS_X\otimes\Q)\cap H_2(X)$ be its primitive hull.
The quotient $\tilde\CS_X/\CS_X$ is equal to the torsion
$\Tors\MW(X)$ of the Mordell--Weil group of~$X$,
see~\cite{Shioda}.

The orthogonal complement $\CT_X:=\CS_X^\perp$ is called the
\emph{\rom(stable\rom) transcendental lattice} of~$X$.
Note that $\CS_X$ is
nondegenerate; hence so is $\CT_X$.

The collection $(H_2(X),\Gs_X,[E])$, considered up to
auto-isometries
of $H_2(X)$ preserving~$[E]$ and $\Gs_X$ as a set,
is called the
\emph{homological type} of~$X$. If $\CS_X$ is primitive,
the homological type is
determined by the combinatorial type of the singular fibers
of~$X$, the lattice~$\CT_X$, and the anti-isometry
$q\:\discr\CS_X\to\discr\CT_X$
defining the extension $H_2(X)\supset\CS_X$,
see~\ref{s.extension}.

\subsection{The skeleton~$\Gamma_X$}\label{s.GammaX}
Let~$X$ be an extremal elliptic surface over a base~$B$.
Define its
\emph{skeleton} as the embedded bipartite graph
$\Gamma_X:=j_X\1[0,1]\subset B$. The
pull-\allowbreak backs
of~$0$ and~$1$ are
called, respectively, \black-- and \white-vertices of~$\Gamma_X$.
(Thus, $\Gamma_X$ is the \emph{dessin d'enfants} of~$j_X$ in the
sense of Grothendieck; however, we reserve the word `dessin' for
the more complicated graphs describing arbitrary, not necessarily
extremal, surfaces, see~\cite{degt.kplets}.)
Since $X$ is extremal, $\Gamma_X$ has the following properties:
\roster
\item\local{Sk.disk}
each region of~$\Gamma_X$ (\ie., component of $B\sminus\Gamma_X$)
is a topological disk;
\item\local{Sk.valency}
the valency of each \black-vertex is~$\le3$, the valency of each
\white-vertex is~$\le2$.
\endroster
In particular, it follows that $\Gamma_X$ is connected.

The skeleton~$\Gamma_X$ determines~$j_X$; hence the pair
$(\Gamma_X,\mh_X)$ determines~$X$. (Here, it is important that $B$
is considered as a \emph{topological} surface; its analytic
structure is given by the Riemann existence theorem.)

\paragraph\label{s.no.E}
From now on, we will speak about extremal surfaces
{\em without $\tE$ type singular fibers}. In this case, all
\black-vertices of~$\Gamma_X$ are of valency~$3$ and all its
\white-vertices are of valency~$2$. Hence, the \white-vertices can
be disregarded (with the convention that a \white-vertex is to be
understood at the center of each edge connecting two
\black-vertices). Furthermore, in view of
condition~\itemref{s.GammaX}{Sk.disk} above, one can also
disregard the underlying surface~$B$ and retain the ribbon graph
structure of~$\Gamma_X$ only. For future references, we restate
the definition:
\roster
\item"$(*)$"
$\Gamma_X$ is a ribbon graph with all vertices of valency~$3$.
\endroster
Under the assumptions, the surface~$B$
containing~$\Gamma$ is reconstructed from the
ribbon graph structure. Its genus is called the \emph{genus}
of~$\Gamma$.

In Subsection~\ref{s.mhX} below, we explain that the homological
invariant~$\mh_X$ can be described in terms of an orientation
of~$\Gamma_X$, reducing an extremal elliptic surface to an
oriented $3$-regular ribbon graph.\mnote{former 2.3.2 removed;
a short remark and a ref here}

\section{The topological aspects\label{S.topology}}

\subsection{The notation}\label{s.notation}
Consider a Jacobian elliptic surface $\pr\:X\to B$ over a base~$B$
of genus~$g$. Let $E\subset X$ be the section of~$X$, and denote
by $F_1,\ldots,F_r$ its singular fibers. Let $S=\bigcup_iF_i$,
$i=1,\ldots,r$.

Recall that \emph{stable} are the singular fibers of~$X$ of
type~$\tA_0^*$ or~$\tA_p$, $p\ge1$. One has
$H_1(F_i)=\Z$ is $F_i$ is stable and $H_1(F_i)=0$ otherwise.

For each $i=1,\ldots,r$,
pick a regular neighborhood $N_i$ of~$F_i$ of the
form $\pr\1U_i$, where $U_i\subset B$ is a small disk about
$\pr F_i$. Let $N_S=\bigcup N_i$, $i=1,\ldots,r$. Let, further,
$N_E$ be a tubular neighborhood of~$E$. We assume
$N_E$ and all~$N_i$
so small that $N:=N_E\cup N_S$ is a regular
neighborhood of $E\cup S$. Thus, the spaces~$N$, $N_E$, and~$N_i$
contain, respectively, $E\cup S$, $E$, and~$F_i$ as strict
deformation retracts.

Denote by~$\Xcirc$ the closure of $X\sminus N$ and decompose the
boundary $\dX$ into the union
$\partial_E\Xcirc\cup\partial_S\Xcirc$,
$\partial_S\Xcirc:=\bigcup\partial_i\Xcirc$, where
$\partial_E\Xcirc:=\dX\cap N_E$ and
$\partial_i\Xcirc:=\dX\cap N_i$, $i=1,\ldots,r$.
Since $\dX=\partial N$, we will use the same notation
$\partial_\bullet N$ for the corresponding parts of the boundary
of~$N$, so that $\partial_\bullet N=\partial_\bullet\Xcirc$.

We also use the notation $\CS_X$, $\tilde\CS_X$, and~$\CT_X$
introduced in~\ref{s.hom.type}.

\subsection{Tubular neighborhoods}
First,
recall that the inclusion $E\into X$ induces isomorphisms, see,
\eg.,~\cite{FM},
$$
H_1(E)@>\cong>>H_1(X),\qquad
H^1(X)@>\cong>>H^1(E).
\eqtag\label{eq.E=X}
$$
The inverse isomorphisms are induced by the projection
$\pr\:X\to B$ and
the obvious
identification $E=B$.

Consider a singular fiber~$F_i$, $i=1,\ldots,r$. The boundary
$\partial_iN=\partial N_i\sminus\operatorname{interior}N$\mnote{
`$\operatorname{interior}$'}
is fibered over the circle
$\partial U_i$, the fiber being a punctured torus~$\Fcirc$.
Denote by~$\mm_i$ and $\mm^*_i$ the monodromy of this fibration in
$H_1(\Fcirc)$ and $H^1(\Fcirc)$, respectively. One has
$$
\gather
H_2(\partial_iN)=\Ker[(\mm_i-\id)\:H_1(\Fcirc)\to H_1(\Fcirc)],
\eqtag\label{eq.dX.h}\\
H^2(\partial_iN)=\Coker[(\mm^*_i-\id)\:H^1(\Fcirc)\to H^1(\Fcirc)].
\eqtag\label{eq.dX.ch}
\endgather
$$
All monodromies~$\mm_i$ are known, see, \eg.,~\cite{FM} or
Example~\ref{ex.region} below.
In particular, $\mm_i$ has invariant vectors if and only if $F_i$
is a stable singular fiber.
Thus, $H_2(\partial_SN)$ is a free group
and one has
$$
\rank H_2(\partial_SN)=
\rank H_1(S)=
 \text{number of stable singular fibers of~$X$}.
\eqtag\label{eq.H1dX}
$$

\paragraph\label{s.YdY}
Let~$Y$ be an oriented
$4$-manifold with boundary. Recall that, if $H_1(Y)$
is torsion free (or, equivalently, $H^2(Y)$ is torsion free),
then $H_2(Y,\partial Y)=H^2(Y)=(H_2(Y))^*$ and the relativization
homomorphism $\rel\:H_2(Y)\to H_2(Y,\partial Y)$
coincides with the
homomorphism associated with the
intersection index form, see Subsection~\ref{s.lattices}.
In particular, one has isomorphisms
$\Tors\Coker\rel=\Tors H_1(\partial Y)=\discr H_2(Y)$.
(The resulting $(\Q/\Z)$-valued bilinear form on
$\Tors H_1(\partial Y)$ is called the \emph{linking coefficient
form}; it can be defined geometrically in terms of $\partial Y$
only.)

Since $H_1(N)=H_1(S\cup E)$ is torsion free and
$H_2(N)=\CS_X/\Rad$, one has
$$
\discr\CS_X=\Tors H_1(\partial N)=\Tors H^2(\partial N).
$$

\lemma\label{d=dS}
The inclusion homomorphism
$H^2(\partial N)\to H^2(\partial_SN)$
restricts to an isomorphism
$\Tors H^2(\partial N)=\Tors H^2(\partial_SN)$.
\endlemma

\proof
Denote $\partial'N_S=\partial N_S\cap N_E$ and consider the
commutative diagram
$$
\CD
H_2(N)@>\rel_1>>H_2(N,\partial N)
 @>\partial_1>>H_1(\partial N)\\
@VVV@VVV@VVV\\
H_2(N_S,\partial'N_S)@>\rel_2>>H_2(N_S,\partial N_S)
 @>\partial_2>>H_1(\partial N_S,\partial'N_S),
\endCD
$$
where the rows are fragments of exact sequences of pairs and
vertical arrows are induced by appropriate inclusions, the
rightmost arrow being Poincar\'e dual to the homomorphism in
question. The cokernels $\Coker\partial_i$, $i=1,2$, belong to the
free groups $H_1(N)$ and $H_1(N_S,\partial'N_S)$, respectively;
hence all torsion elements come from the cokernels $\Coker\rel_i$.
It remains to observe that
$$
\CS_X=H_2(N)/\Rad=\CU\oplus(H_2(N_S)/\Rad),
$$
hence
$\Tors\Coker\rel_1=\discr H_2(N)=\discr H_2(N_S)=\Coker\rel_2$.
To establish the last equality, notice that, for each singular
fiber~$F_i$, there is a decomposition (not orthogonal)
$H_2(N_i,\partial'N_i)=H_2(N_i)\oplus\Z[E_i,\partial E_i]$,\mnote{$N_i$}
where
$E_i=E\cap N_i$; hence one can identify $H_2(N_i,\partial'N_i)$
with $(H_2(N_i)/\Rad)\oplus\CU$.
\endproof

The advantage of Lemma~\ref{d=dS} is the fact that the
isomorphisms $\discr H_2(N_i)=\Tors H^2(\partial_iN)$ are local:
they can be computed in terms of the
topological types of the singular fibers
of~$X$.

\subsection{The homology of $\Xcirc$}
In this subsection, we compute the invariants~$\CT_X$ and
$\Tors\MW(X)$ of an arbitrary Jacobian elliptic surface~$X$ in
terms of the \hbox{(co-)}\allowbreak homology of~$\Xcirc$.

\lemma\label{lem.HX}
The group $H_2(\Xcirc)$ is free and there is a short exact
sequence
$$
0@>>>\Rad H_2(\Xcirc)@>>>H_2(\Xcirc)@>>>\CT_X@>>>0,
$$
so that $\CT_X=H_2(\Xcirc)/\Rad$.
Furthermore,
the homomorphism $H_2(\partial_S\Xcirc)\to H_2(\Xcirc)$ induced by
the inclusion establishes an isomorphism
$H_2(\partial_S\Xcirc)=\Rad H_2(\Xcirc)$.
\endlemma

\proof
The first statement is an immediate consequence
from the Poincar\'e duality
$H_2(\Xcirc)=H^2(\Xcirc,\dX)$
and the exact sequence
$$
H^1(X)@>>>H^1(N)@>\partial>>H^2(\Xcirc,\dX)@>>>H^2(X)@>>>H^2(N);
$$
the kernel of the last homomorphism is
$\CT_X\subset H^2(X)=H_2(X)$,
and
the cokernel $H^1(N)/H^1(X)=H^1(S)$ is free, \cf.~\eqref{eq.E=X}.
As another consequence,
the rank
$\rank\Im\partial$ equals the number of
stable singular fibers of~$X$.

The homomorphism~$\partial$ above is Poincer\'e dual to~$\partial$
in the following commutative diagram:
$$
\CD
H_3(X,\Xcirc)@>\partial>>H_2(\Xcirc)\\
@|@AA\inj_*A\\
H_3(N,\dX)@>>>H_2(\dX).
\endCD
$$
It follows that
$\Im\partial\subset\Im\inj_*\subset\Rad H_2(\Xcirc)$. (Classes
coming from the boundary are always in the kernel of the
intersection index form.)
Since
$\CT_X$ is nondegenerate, both inclusions are equalities.

Finally, consider the exact sequence
$$
H_2(\partial_S\Xcirc)@>>>H_2(\dX)@>>>
 H_2(\dX,\partial_S\Xcirc)@>\partial>>H_1(\partial_S\Xcirc).
$$
One has
$H_*(\dX,\partial_S\Xcirc)=H_*(E',\partial E')\otimes H_*(S^1)$,
where $E'=E\sminus N_S$, and it is easy to see that
$\Ker\partial=H_1(E',\partial E')\otimes H_1(S^1)$ and that each
element of this kernel lifts to a class in $H_2(\dX)$ that
vanishes in $H_2(\Xcirc)$. (If
$\Ga$ is a relative $1$-cycle in $(E',\partial E')$,
the lift is the boundary of $\pr\1\pr\Ga\sminus N_E$.) Thus,
the image of $H_2(\dX)$ in $H_2(\Xcirc)$ coincides with that of
$H_2(\partial_S\Xcirc)$. Since the ranks of
$H_2(\partial_S\Xcirc)$ and its image coincide (both equal to the
number of stable singular fibers of~$X$), the inclusion induces an
isomorphism.
\endproof

\lemma\label{lem.H*X}
There is an exact sequence
$$
0@>>>\CS_X\to H_2(X)\to H^2(\Xcirc)\to H_1(S)\to0.
$$
In particular,
$\Tors H^2(\Xcirc)=\tilde\CS_X/\CS_X=\Tors\MW(X)$.
\endlemma

\proof
The statement follows from the Poincar\'e duality
$H^2(\Xcirc)=H_2(\Xcirc,\dX)$, the exact sequence
$$
H_2(N)@>>>H_2(X)@>>>H_2(\Xcirc,\dX)@>>>
 H_1(N)@>>>H_1(X),
$$
and the fact that
$\Ker[H_1(N)\to H_1(X)]=H_1(S)$, \cf.~\eqref{eq.E=X}.
\endproof

Assume that $\CS_X$ is primitive in $H_2(X)$, \ie.,
$\tilde\CS_X=\CS_X$. Then, due to~\ref{s.YdY} and
Lemma~\ref{lem.H*X}, there is an isomorphism
$\discr\CT_X=\Tors H^2(\dX)$, which gives rise to an isomorphism
$\discr\CT_X=\Tors H^2(\partial_S\Xcirc)$, see Lemma~\ref{d=dS}.

\lemma\label{lem.q}
If $\CS_X$ is primitive in $H_2(X)$, then the anti-isometry
$q\:\discr\CT_X\to\discr\CS_X$ defining the homological type
of~$X$, see~\ref{s.hom.type},
can be identified with
the composition $j\1\circ i$ of the isomorphisms
$$
\discr\CT_X@>i>\cong>\Tors H^2(\partial_S\Xcirc)
 @<j<\cong<\discr\CS_X
$$
induced by the inclusions $\partial_S\Xcirc\into\Xcirc$ and
$\partial_S\Xcirc\into N_S$.
\endlemma

\proof
Using Lemma~\ref{d=dS}, one can replace $\partial_S\Xcirc$
with~$\dX$. Then the statement follows from the Mayer--Vietoris
exact sequence
$$
H^2(X)@>>>H^2(N)\oplus H^2(\Xcirc)@>>>H^2(\dX)
$$
and the definition of~$q$.
\endproof

\subsection{The counts}\label{s.counts}
We conclude this section with a few counts.

\paragraph\label{s.Sigma}
Let~$X$ be an extremal
elliptic surface
over a curve~$B$ of genus~$g$, and
let $\Gamma=\Gamma_X\subset B$ be the skeleton
of~$X$.
Assume that all singular fibers of~$X$ are of type~$\tA_0^*$,
$\tA_p$, $p\ge1$, or~$\tD_q$, $q\ge5$, and
denote by~$t$ the number of $\tD$~type fibers.
Let $\chi(X)=6(k+t)$. (Recall that $12\mathrel|\chi(X)$.)
Then
the quotient $X/\pm1$ blows down to a ruled surface~$\Sigma$
over~$B$ with an exceptional section~$E$ with $E^2=-(k+t)$. The
ramification locus of the projection $X\to\Sigma$ is the union
$C\cup E$, where~$C$ is a certain
\emph{trigonal curve} (\ie., a curve
disjoint from~$E$ and intersecting each generic fiber of the
ruling at three points) with simple singularities only.

The surface~$X$ is diffeomorphic to the double covering
$X'\to\Sigma$ ramified at~$E$ and a nonsingular trigonal
curve~$C'$. Using this fact and taking into
account~\eqref{eq.E=X}, one can easily compute the inertia indices
$\Gs_\pm$ of the intersection index form on~$H_2(X)$:
$$
\Gs_+(X)=k+t+2g-1,\qquad
\Gs_-(X)=5k+5t+2g-1.
$$

\paragraph\label{s.v-e-r}
Let $\Gamma=\Gamma_X\subset B$ be the skeleton of~$X$.
The numbers of vertices, edges, and regions of~$\Gamma$ are,
respectively,
$$
v=2k,\qquad e=3k,\qquad r=k+2-2g.
$$
The latter count~$r$ is also the number of singular fibers
of~$X$.
The `total Milnor number' of the singular
fibers of~$X$ is given by
$\mu=2g+5k+5t-2$. (Indeed, each $n$-gonal region~$R$ contributes
$(n-1)$ or $(n+4)$ depending on whether $R$ contains an $\tA$
or~$\tD$ type fiber. The total number of corners of the regions
is~$6k$.) Taking into account Lemma~\ref{lem.HX}
and~\eqref{eq.H1dX}, one arrives at the following statement.

\lemma\label{counts}
In the notation above, $\CT_X$ is a positive definite lattice of
rank $k+t+2g-2$.
Furthermore, one has $\rank\Rad H_2(\Xcirc)=k-t+2-2g$, and
$H_2(\Xcirc)$ is a positive semi-definite lattice of rank~$2k$.
\qed
\endlemma

\Remark\label{rem.positive}
The assertion that
the lattice~$\CT_X$ is positive definite still holds if $X$
has type~$\tE$ singular fibers. In fact,
this property can be taken for the
definition of an extremal elliptic surface.
\endRemark

\section{The main theorem\label{S.main}}

\subsection{Skeletons}
To ease the further exposition, we redefine a skeleton
in the sense of~\ref{s.no.E}$(*)$
as a set of ends of its edges. However, we will make no
distinction between a skeleton in the sense of
Definition~\ref{def.skeleton} below and its geometric
realization.

\definition\label{def.skeleton}
A \emph{skeleton} is a collection $\Gamma=(\CE,\op,\nx)$, where
$\CE$ is a finite set, $\op\:\CE\to\CE$ is a free involution, and
$\nx\:\CE\to\CE$ is a free automorphism of order~$3$.
The orbits of~$\op$ are called
the\mnote{one `the' removed} \emph{edges} of~$\Gamma$,
and the orbits of~$\nx$ are called its \emph{vertices}.
(Informally, $\op$ assigns to an end the other end\mnote{`end'}
of the same
edge, and $\nx$ assigns the next end at the same vertex with
respect to its cyclic order.)
\enddefinition

\paragraph\label{s.orientation}
According to this definition, the sets of edges and vertices
of a skeleton~$\Gamma$ can be referred to as $\CE\!/\!\op$ and $\CE\!/\!\nx$,
respectively. An \emph{orientation} of~$\Gamma$ is a section
${+}\:\CE\!/\!\op\to\CE$ of~$\op$, sending each edge~$e$ to its
\emph{head}~$e^+$.
Given such a section, its composition with~$\op$ sends each edge~$e$
to its \emph{tail}~$e^-$.
It is worth mentioning that, from this point of view,
a \emph{marking} of~$\Gamma$ in the sense of~\cite{degt.kplets}
is merely
a section $\m1\:\CE\!/\!\nx\to\CE$ of~$\nx$,
sending each vertex to the
first edge end attached to it.
Then the sections $\m2:={\nx}\circ\m1$
and $\m3:={\nx^2}\circ\m1$ send a vertex to the second
and third edge ends, respectively.

The elements $\op$ and~$\nx$ of order~$2$ and~$3$, respectively,
generate the modular group $\PSL(2,\Z)\cong\CG2*\CG3$, which acts
on~$\CE$. A skeleton is \emph{connected} if this action is
transitive. Recall that each element $w\in\PSL(2,\Z)$ can be
uniquely represented by a reduced word $w_1w_2w_3\ldots$
of the form
$\op\nx^{\pm1}\op\ldots$
or $\nx^{\pm1}\op\nx^{\pm1}\ldots$. The length of this
word is called the \emph{length} of~$w$.

\Remark\label{rem.computer}
It is worth mentioning that Definition~\ref{def.skeleton} results
in a computer friendly presentation of~$\Gamma$: it is given by
two permutations~$\op$ and~$\nx$, the former splitting into a
product of cycles of length~$2$, the latter, into a product of
cycles of length~$3$. Certainly, this description
is equivalent to the
presentation of the ramified covering $B\to\Cp1$ defined
by~$\Gamma$ by its Hurwitz system.
\endRemark

\definition\label{def.path}
A \emph{path} in a skeleton $\Gamma=(\CE,\op,\nx)$
can be defined as a
pair $\Gg=(\Ga,w)$, where $\Ga\in\CE$ and $w\in\PSL(2,\Z)$.
If $w$ is a positive power of $\nx\1\op$, then $\Gg$ is
called a \emph{left turn path}
(\cf. Figure~\ref{fig.paths}, left, in Subsection~\ref{s.all.A}
below).
The \emph{endpoint} of~$\Gg$ is the element $w(\Ga)\in\CE$.
If the length of~$w$ is even and $w(\Ga)=\Ga$,
the path is called a \emph{loop}.
\enddefinition

\paragraph\label{s.path=chain}
Representing~$w$ by a reduced word $w_r\ldots w_1$, one can
identify a path $(\Ga,w)$ with a sequence $(\Ga_0,\ldots,\Ga_r)$,
where $\Ga_0=\Ga$ and $\Ga_i=w_i(\Ga_{i-1})$ for $i\ge1$.

\paragraph\label{s.region}
A \emph{region} of a skeleton~$\Gamma$ can be defined as an
orbit of the cyclic subgroup of
$\PSL(2,\Z)$ generated by $\nx\1\op$. Given an $n$-gonal region~$R$,
$n\ge1$, and an
element $\Ga_0\in R$, the \emph{boundary} $\dR$ is the left turn
path of length~$2n$ starting at~$\Ga_0$. It is a loop. In the
sequence $(\Ga_0,\Ga_1,\ldots,\Ga_{2n}=\Ga_0)$ representing~$\dR$,
each even term $\Ga_{2i}$ is an element of~$R$, and each odd term
has the form $\Ga_{2i+1}=\op\Ga_{2i}$.

Patching the boundary of each region of~$\Gamma$ with a disk, one
obtains the surface~$B$ containing~$\Gamma$. Hence, the genus
$g(\Gamma)$
of~$\Gamma$, see~\ref{s.no.E}, is given by
$$
2-2g(\Gamma)=\mathop\#(\CE\!/\!\nx)
 -\mathop\#(\CE\!/\!\op)+\mathop\#(\CE\!/\!\nx\1\op).
$$

\subsection{The homological invariant}\label{s.mhX}
Let\mnote{part of \ref{s.tripods} moved here}
$\CH=\Z a\oplus\Z b$ with the skew-symmetric bilinear form
$\bigwedge^2\CH\to\Z$ given by $a\cdot b=1$. Introduce the
isometries $\X,\Y\:\CH\to\CH$ given (in the standard basis $\{a,b\}$) by
the matrices
$$
\X=\bmatrix-1&1\\-1&0\endbmatrix,\qquad
\Y=\bmatrix0&-1\\1&\phantom{-}0\endbmatrix.
$$
One has $\X^3=\id$ and $\Y^2=-\id$. If $c=-a-b\in\CH$, then $\X$ acts
\via\
$$
(a,b)\overset \X\to\longmapsto(c,a)
 \overset \X\to\longmapsto(b,c)
 \overset \X\to\longmapsto(a,b).
$$
It is well known that $\X$ and~$\Y$ generate the group $\SL(2,\Z)$
of isometries of~$\CH$. We fix the notation $\CH$, $a$, $b$, $c$
and $\X$, $\Y$
throughout the paper.\mnote{notation fixed}

Let\mnote{former 4.2.8 moved here and extended}
$\pr\:X\to B$ be an elliptic surface with singular fibers of
type~$\tA{_0^*}$, $\tA{_p}$, $p\ge1$, or $\tD{_q}$, $q\ge5$, only.
We use the results of~\cite{degt.kplets} to describe the
homological invariant of~$X$
in terms of the skeleton $\Gamma=\Gamma_X$.
More precisely, we describe the
monodromy in $H_1(\text{fiber})$ of the locally trivial fibration
$\pr\:\pr\1\Gamma\to\Gamma$.

Consider the double covering $X\to\Sigma$ ramified at $C\cup E$,
see~\ref{s.Sigma}. Pick a vertex $v$ of~$\Gamma$, let $F_v$ be the
fiber of~$X$ over~$v$, and let $\barF_v$ be its projection
to~$\Sigma$. Then, $F_v$ is the double covering of~$\barF_v$
ramified at $\barF_v\cap(C\cup E)$ (the three black points in
Figure~\ref{fig.basis} and~$\infty$).

\midinsert
\centerline{\picture{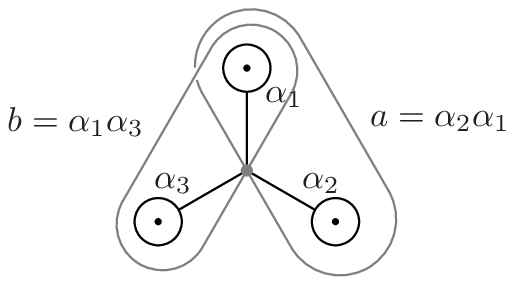}}
\figure
The basis in $H_1(F_v)$
\endfigure\label{fig.basis}
\endinsert

In the presence of a trigonal curve, $\Sigma$ has a well defined
zero section (the fiberwise barycenter of the points of the
curve with respect to the canonical $\C^1$-affine structure in the
open fibers $\barF\sminus E$).
Let $\bz_v\in\barF_v$ be the value of the zero section at a vertex~$v$
of~$\Gamma$. For each vertex~$v$, pick and fix one of the two
pull-backs of~$\bz_v$ in~$F_v$; denote it by~$z_v$.
The collection $\{z_v\}$, $v\in\CE\!/\!\nx$, is called a
\emph{reference set}.

Choose a marking at~$v$ and let $\{\Ga_1,\Ga_2,\Ga_3\}$ be the
canonical basis for
the group $\pi_1(\barF_v\sminus(C\cup E),\bz_v)$ defined by
this marking (see~\cite{degt.kplets} and Figure~\ref{fig.basis};
unlike~\cite{degt.kplets}, we take~$\bz_v$ for the reference point;
this choice removes the ambiguity in the
definition of canonical basis).
Then $H_1(F_v)=\pi_1(F_v,z_v)$ is generated by the
lifts $a_v=\Ga_2\Ga_1$ and $b_v=\Ga_1\Ga_3$ (the two grey cycles in
the figure), and one can use
the map $a_v\mapsto a$, $b_v\mapsto b$
to identify
$H_1(F_v)$ with~$\CH$.

In the sequel, we consider a separate copy~$F_\Ga$ of~$F_v$ for
each edge end $\Ga\in v$.

\definition\label{def.id}
The \emph{canonical identification} is the isomorphism
$H_1(F_\Ga)\to\CH$ constructed above using the marking at~$v$
defined \via\ $\Ga=\m1(v)$.
\enddefinition

\lemma\label{lem.X}
Under the canonical identification, the
identity map $F_\Ga\to F_{\nx\Ga}$, regarded as an automorphism
of~$\CH$, is given by~$\X\1$.
\endlemma

\proof
This map is the change of basis from $\{a,b\}$
to $\{c,a\}$.
\endproof

\lemma\label{lem.Y}
Let~$u$ and~$v$ be two vertices \rom(not necessarily distinct\rom)
connected by an edge~$e$, and let $\Ga\in u$ and $\Gb\in v$ be the
respective ends of~$e$.
Under the canonical identifications over~$u$ and~$v$,
the monodromy $H_1(F_\Ga)\to H_1(F_\Gb)$
along~$e$, regarded as an automorphism of~$\CH$, is given by
$\pm\Y$.
\endlemma

\proof
This monodromy
is a lift of monodromy~$m_{1,1}$
in~\cite{degt.kplets}; geometrically (in~$\Sigma$), the black
ramification point
surrounded by~$\Ga_1$ crosses the segment connecting
the ramification points surrounded by~$\Ga_2$
and~$\Ga_3$.
\endproof

The sign $\pm1$ in Lemma~\ref{lem.Y}
depends on the homological invariant~$\mh_X$
\emph{and} on the choice of a
reference set. The monodromy from~$v$
to~$u$ is $(\pm\Y)\1=\mp\Y$.

\definition\label{def.orientation}
Given an elliptic surface~$X$ as above and a reference set
$\{z_v\}$, $v\in\CE\!/\!\nx$,
we define an orientation of~$\Gamma$ as follows:
an edge~$e$ is oriented so
that the monodromy $H_1(F_{e^-})\to H_1(F_{e^+})$ along~$e$ be
given by $+\Y$.
\enddefinition

Changing the lift~$z_v$ over a vertex~$v$ to the other one
results in a change of sign of the canonical identification
$H_1(F_\Ga)\to\CH$ for each end $\Ga\in v$.
As a consequence, each monodromy
starting or ending at~$v$ changes sign.
Thus,
two orientations of~$\Gamma$
give rise to the same monodromy over~$\Gamma$ if and
only if
they are obtained from each other
by the following
operation: pick a subset $V$ of the set of vertices of~$\Gamma$
and reverse the orientation of each edge that has exactly
one end in~$V$. Summarizing, one arrives at the following
statement.

\lemma\label{lem.orientation}
An extremal elliptic surface~$X$
without $\tE$ type singular fibers is determined up to isomorphism
by an oriented
ribbon
graph~$\Gamma_X$
as in~\ref{s.no.E}$(*)$. Conversely, oriented
ribbon graph~$\Gamma_X$ is
determined by~$X$
up to isomorphism and a change of orientation just described.
\endlemma

\proof
If $X$ is extremal and without $\tE$ type fibers, then $\Gamma_X$
is a strict deformation retract of~$\Bcirc$ and the monodromy
over~$\Gamma_X$ determines~$\mh_X$.
\endproof

\subsection{The tripod calculus}\label{s.tripods}
Let $\Gamma=(\CE,\op,\nx,+)$ be a connected oriented skeleton.
Place a copy~$\CH_\Ga$
of~$\CH$ at each element $\Ga\in\CE$, and let
$\CH\gtimes\Gamma=\bigoplus\CH_\Ga$, $\Ga\in\CE$.
For a vector $h\in\CH\gtimes\Gamma$, we denote by~$h_\Ga$ its
projection to~$\CH_\Ga$, $\Ga\in\CE$; for a vector $u\in\CH$ and
element $\Ga\in\CE$, denote by $u\gtimes\Ga\in\CH\gtimes\Gamma$
the vector whose only nontrivial projection is
$(u\gtimes\Ga)_\Ga=u$.\mnote{`$u$'}
Convert $\CH\gtimes\Gamma$
to a rational lattice by letting
$$
h^2=-\frac13\sum_{\Ga\in\CE}h_\Ga\cdot\X h_{\nx\Ga},\quad
h\in\CH\gtimes\Gamma,
\eqtag\label{eq.form}
$$
where $\,\cdot\,$ stands for the product in~$\CH$. Let
$\CH_\Gamma$ be the sublattice of $\CH\gtimes\Gamma$ subject to
the following relations:
\roster
\item\local{d.vertex}
$h_{\Ga}+\X h_{\nx\Ga}+\X^2h_{\nx^2\Ga}=0$ for each element
$\Ga\in\CE$;
\item\local{d.edge}
$h_{e^+}+\Y h_{e^-}=0$ for each edge $e\in\CE\!/\!\op$.
\endroster
Similarly, consider the dual group
$\CH^*\gtimes\Gamma=\bigoplus\CH^*_\Ga$, $\Ga\in\CE$, where
$\CH^*_\Ga$ is a copy of
the dual group~$\CH^*$, and define $\CH^*_\Gamma$ as the
quotient of $\CH^*\gtimes\Gamma$ by the subgroup spanned by the
vectors of the form
\roster
\item[3]\local{d*.vertex}
$u\gtimes\Ga+\X^*u\gtimes(\nx\Ga)+(\X^*)^2u\gtimes(\nx^2\Ga)$
for each $u\in\CH^*$ and $\Ga\in\CE$;
\item\local{d*.edge}
$u\gtimes e^++\Y^*u\gtimes e^-$ for each
$u\in\CH^*$ and $e\in\CE\!/\!\op$.
\endroster
(Here, $\X^*,\Y^*\:\CH^*\to\CH^*$ are the adjoint of $\X$, $\Y$.)
It is easy to see that $\CH_\Gamma$ annihilates the subgroup
spanned by~\loccit{d*.vertex}, \loccit{d*.edge}, inducing a
pairing $\CH_\Gamma\otimes\CH^*_\Gamma\to\Z$.
(Observe that the maps $h\mapsto h_\Ga\in\CH$ and
$u\mapsto u\gtimes\Ga\in\CH^*\gtimes\Gamma$ are adjoint to each
other.)
Note that, in general, $\CH^*_\Gamma\ne(\CH_\Gamma)^*$, as
$\CH^*_\Gamma$ may have torsion.

\Remark\label{rem.1}
Since $\X^3=\id$, in relation~\loccit{d.vertex} above
it suffices to pick a
marking $\m1\:\CE\!/\!\nx\to\CE$, see~\ref{s.orientation},
and consider
one relation
\roster
\item[5]\local{d.vertex.alt}
$h_{\m1(v)}+\X h_{\m2(v)}+\X^2h_{\m3(v)}=0$ for each vertex
$v\in\CE\!/\!\nx$.
\endroster
Furthermore, since $\X$ is an isometry,
the restriction to~$\CH_\Gamma$ of the
quadratic form given by~\eqref{eq.form} can be simplified to
$$
h^2=-\sum_{v\in\CE\!/\!\nx}h_{\m1(v)}\cdot\X h_{\m2(v)},
\qquad h\in\CH\gtimes\Gamma.
\eqtag\label{eq.form.alt}
$$
This expression (when restricted to $\CH_\Gamma$) does not depend
on the marking.
\endRemark

Now, let $X=X_\Gamma$ be the extremal elliptic surface
defined by~$\Gamma$, see~Lemma~\ref{lem.orientation}.
Next theorem
computes the (co-)homology of $\Xcirc_\Gamma$,
see~\ref{s.notation}, in terms of~$\Gamma$.

\theorem\label{tripods}
There are isomorphisms $H_2(\Xcirc_\Gamma)=\CH_\Gamma$ and
$H^2(\Xcirc_\Gamma)=\CH^*_\Gamma$.
The former takes the intersection index form to the form given
by~\eqref{eq.form}\rom; the latter takes the Kronecker product to
the pairing $\CH_\Gamma\otimes\CH^*_\Gamma$ defined above.
\endtheorem

\proof
Replace $\Xcirc$ with its strict deformation retract
$X':=\pr\1\Gamma\sminus N_E$; it fibers over~$\Gamma$ with the
fiber punctured torus. Subdivide~$\Gamma$ into cells by taking its
\black-- and \white-vertices for $0$-cells and half edges
(\ie., edges of the form $\edge$)
for
$1$-cells, and let~$X'_0$ be the pull-back of the $0$-skeleton
of~$\Gamma$. Then, in the exact sequence
$$
H_2(X'_0)@>>>H_2(X')@>>>H_2(X',X'_0)@>\partial>>H_1(X'_0)
$$
of pair $(X',X'_0)$ one has $H_2(X'_0)=0$; hence
$H_2(\Xcirc)=H_2(X')=\Ker\partial$.\mnote{the rest of the proof expanded;
references to \ref{s.mhX}}

Pick a marking of~$\Gamma$, see~\ref{s.orientation}, and a
reference set $\{z_\Ga\}$, $\Ga\in\CE\!/\!\nx$,
with respect to which $\mh_X$ defines the given
orientation of~$\Gamma$, see
Definition~\ref{def.orientation}. Note that, for each fiber~$F$,
the inclusion $\Fcirc\into F$ induces an isomorphism
$H_1(\Fcirc)=H_1(F)$.

The half edges of~$\Gamma$ are in a one-to-one correspondence with
the elements of $\CE$, and, under the canonical identifications,
see Definition~\ref{def.id}, the group $H_2(X',X'_0)$ splits into
direct sum
$\bigoplus_{\Ga\in\CE}H_1(F_\Ga)\otimes H_1(I_\Ga,\partial I_\Ga)
 =\bigoplus_{\Ga\in\CE}\CH\otimes\Z$, where $I_\Ga$ is the half
edge containing~$\Ga$. To establish an isomorphism
$H_1(I_\Ga,\partial I_\Ga)=\Z$, we
use the fundamental class $[I_\Ga,\partial I_\Ga]$ corresponding
to the orientation of~$I_\Ga$
towards its
\black-vertex.
In other words, for each $\Ga\in\CE$, we consider a direct summand
$$
\CH_\Ga:=
H_1(F_\Ga)\otimes\Z[\oedge],
\quad H_1(F_\Ga)=\CH.
\eqtag\label{eq.Halpha}
$$
Thus, there is a canonical isomorphism
$H_2(X',X_0')=\CH\gtimes\Gamma$.

For each \black-vertex~$v$, identify $H_1(F_v)$ with~$\CH$ using
the chosen marking, so that $H_1(F_v)=H_1(F_{\m1(v)})$.
Then the composition
$H_2(X',X'_0)\to H_1(X'_0)\to H_1(F_v)=\CH$ of
the boundary operator~$\partial$ and the
projection to $H_1(F_v)$ is given by the left hand side
of~\iref{rem.1}{d.vertex.alt} at~$v$, see Lemma~\ref{lem.X}.

Finally, a \white-vertex~$w$ of~$\Gamma$ is represented by the
edge~$e$ containing this vertex, and we identify $H_1(F_w)$ with
$H_1(F_{e^+})$ (and further with~$\CH$). Then the composition
$H_2(X',X'_0)\to H_1(X'_0)\to H_1(F_w)=\CH$ is given, up to
sign~$(-1)$, by the left hand side of~\iref{s.tripods}{d.edge}
at~$e$, see Lemma~\ref{lem.Y}.

Thus, after appropriate identifications, $\partial$ is a map
$$
\partial\:\CH\gtimes\Gamma\to
 \bigoplus_{v\in\CE/\nx}\CH\oplus\bigoplus_{e\in\CE/\op}\CH,
\eqtag\label{eq.d}
$$
and its components are given by the left hand sides of the
respective
constraints
\iref{rem.1}{d.vertex.alt} and \iref{s.tripods}{d.edge}
defining $\CH_\Gamma$. Hence one has
$H_2(\Xcirc)=\Ker\partial=\CH_\Gamma$.

The proof for the
cohomology is literally the same, and the interpretation of the
Kronecker product is straightforward.

\midinsert
\centerline{\picture{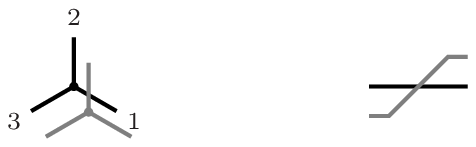}}
\figure
Shift of a marked skeleton
\endfigure\label{fig.shift}
\endinsert

To compute the self-intersection in~$\Xcirc$ of a
$2$-cycle in~$X'$, we
mark~$\Gamma$, shift it in~$\Bcirc$ as shown in Figure~\ref{fig.shift},
left,
and shift the cycle accordingly. Next to each \black-vertex~$v$
of~$\Gamma$, an
intersection point forms; it contributes one term
to~\eqref{eq.form.alt}. (One needs to apply~$\X$ to $h_{\m2(v)}$
in order to bring
$\CH_{\m1(v)}$ and $\CH_{\m2(v)}$
to the same basis, see
Lemma~\ref{lem.X}.)
The shifts do not need to agree, as a
possible intersection point at the middle of an edge of~$\Gamma$,
see Figure~\ref{fig.shift}, right, would not contribute to the
\emph{self}-intersection of a cycle (since self-intersections in
$H_1(\text{fiber})\cong\CH$
are trivial).
\endproof

\corollary\label{independent}
All equations~\iref{s.tripods}{d.edge}
and~\iref{rem.1}{d.vertex.alt} are
linearly independent.
\endcorollary

\proof
This statement follows from Theorem~\ref{tripods} and a
simple dimension count
using Lemma~\ref{counts}.
\endproof

\corollary\label{cor.CT}
There is an isomorphism $\CT_X=\CH_\Gamma/\Rad$.
\endcorollary

\proof
The statement follows from Theorem~\ref{tripods} and
Lemma~\ref{lem.HX}.
\endproof

\corollary\label{cor.MW}
There is an isomorphism $\MW(X_\Gamma)=\Tors\CH^*_\Gamma$.
\endcorollary

\proof
The statement follows from Theorem~\ref{tripods},
Lemma~\ref{lem.HX}, and the fact that the Mordell--Weil group of
an extremal surface has rank~$0$.
\endproof

\Remark
Alternatively,\mnote{new remark}
one can compute $\MW(X_\Gamma)$ in terms
of~$\CH\gtimes\Gamma$ only, \via\
$\MW(X_\Gamma)=\Ext(\Coker\partial,\Z)$, where $\partial$ is
the map given by~\eqref{eq.d}.
\endRemark

\subsection{The monodromy}\label{s.monodromy}
Definition~\ref{def.monodromy} below\mnote{introduction
rewritten; ref to \ref{s.mhX}; signs explained}
is a combinatorial counterpart of the computation of the
homological invariant~$\mh_X$
given by Lemmas~\ref{lem.X} and~\ref{lem.Y}.
Unlike~\ref{s.mhX}, here we are dealing with the groups~$\CH_\Ga$
of $2$-chains, see~\eqref{eq.Halpha}, rather than the groups
$H_1(F_\Ga)$ of $1$-cycles, and we are interested in propagating a
$2$-chain along a path in~$\Gamma$. When following a path, at each
step the orientation in the base is reversed (compared to the
convention $\oedge$ set in~\eqref{eq.Halpha}\,);
it is this fact that explains the
extra sign $-1$ in Definition~\ref{def.monodromy}. In other words,
the sign is chosen so that the parallel transport $\pp{\Gg,h_0}$
defined below be a cycle except over the endpoints of~$\Gg$.
Note that, since loops have even length, the monodromy along a
loop would formally coincide with that given by
Lemmas~\ref{lem.X} and~\ref{lem.Y}.


\definition\label{def.monodromy}
Let $\Gg=(\Ga,w)$ be a path in~$\Gamma$. Represent $w$ by a
reduced word $w_r\ldots w_1$, let $(\Ga_0,\ldots,\Ga_r)$ be the
sequence of vertices of~$\Gg$, and lift $w_i$ and~$w$ to
$\mm_i,\mm=\mm_r\ldots\mm_1\in\SL(2,\Z)$ as follows:
\roster
\item\local{m.vertex}
if $w_i=\nx^{\pm1}$, let $\mm_i=-\X^{\mp1}$,
\item\local{m.edge}
if $w_i=\op$, hence $[\Ga_{i-1},\Ga_i]$ is an edge~$e$
and $\Ga_i=e^\pm$, let $\mm_i=-\Y^{\pm1}$.
\endroster
The map $\mm=\mm_\Gg\:\CH_{\Ga_0}\to\CH_{\Ga_r}$ is called the
\emph{monodromy} along~$\Gg$.
Given a vector $h_0\in\CH$, we define the \emph{parallel transport}
$\pp{\Gg,h_0}\in\CH\gtimes\Gamma$ to be
$\sum_ih_i\gtimes\Ga_i$, where $h_i=\mm_i(h_{i-1})$,
$i=1,\ldots,r$.
\enddefinition

\example\label{ex.region}
The monodromy along the boundary of an $n$-gonal region~$R$
of~$\Gamma$, see~\ref{s.region}, is
$$
\pm(\X\Y)^n=\pm\bmatrix1&n\\0&1\endbmatrix.
$$
Thus, the orientation of~$\Gamma$
determines its type specification
in a simple way: the fiber inside~$R$ is of type~$\tA$ or~$\tD$ if
the sign above is $+$ or~$-$, respectively.
\endexample

\paragraph\label{s.fc}
Let $\Gg=(\Ga,w)$ be a loop, and assume that
the monodromy $\mm_\Gg$ has an invariant vector $h\in\CH_\Ga$.
Then the \emph{fundamental cycle}
$\fc{\Gg,h}:=\pp{\Gg,h}-h\gtimes\Ga$
is an element of~$\CH_\Gamma$.

\example\label{ex.A}
If $R$ is an $n$-gonal region of~$\Gamma$, see~\ref{s.region},
containing an
$\tA$ type singular fiber, then $a$ is invariant under the
monodromy $\mm_{\dR}$, see Example~\ref{ex.region};
hence $\fc{\dR,a}$ is a well defined element of
$\CH_\Gamma=H_2(\Xcirc_\Gamma)$.
(Up to sign, this element does not depend on the choice of the
initial point of~$\dR$.)
Shifting the cycle realizing this element
inside~$R$, one can see that $\fc{\dR,a}\in\Rad H_2(\Xcirc_\Gamma)$.
\endexample

\proposition\label{kernel}
Let $R_1,\ldots,R_{f-t}$ be the regions of~$\Gamma$ containing
its
stable
singular fibers. Then the elements $[\dR_i,a]$,
$i=1,\ldots,f-t$, see Example~\ref{ex.A},
form a basis for the kernel $\Rad\CH_\Gamma$.
\endproposition

\proof
Due to~\eqref{eq.dX.h} and Example~\ref{ex.region}, the elements
in question form a basis for $H_2(\partial_S\Xcirc_\Gamma)$, and the
statement follows from
Lemma~\ref{lem.HX} and Theorem~\ref{tripods}.
\endproof

\paragraph\label{s.cohomology}
Let~$R$ be an $n$-gonal
region of~$\Gamma$. Represent the boundary path $\dR$ by a sequence
$(\Ga_0,\Ga_1,\ldots,\Ga_{n-1})$, see~\ref{s.path=chain},
omitting $\Ga_n=\Ga_0$. Let $\CH\gtimes\dR=\bigoplus\CH^*_i$ be
the direct sum of $n$ copies of~$\CH^*$, one copy for each
vertex~$\Ga_i$, and define the restriction homomorphism
$\res\:\CH^*\gtimes\Gamma\to\CH^*\gtimes\dR$ \via\
$u\gtimes\Ga\mapsto\sum u\gtimes\Ga_i$, the summation running over
all vertices~$\Ga_i$ that are equal to~$\Ga$. (Note that the chain
representing $\dR$ may have repetitions.)

Let $\mm_i^*\:\CH^*_i\to\CH^*_{i-1}$ be the map
adjoint to~$\mm_i$, see Definition~\ref{def.monodromy}.
For~$\mm_n$, we identify $\CH^*_n$ with $\CH^*_0$.
The following statement is straightforward, \cf. the proof of
Theorem~\ref{tripods}; if $\CS_X$ is primitive in $H_2(X_\Gamma)$,
it describes the lattice extension $H_2(X_\Gamma)\supset\CS_X$,
\cf. Lemma~\ref{lem.q}.

\proposition\label{cohomology}
Let~$R$ be an $n$-gonal region of~$\Gamma$ containing a singular
fiber $F_j$ of~$X_\Gamma$. Then there is an isomorphism
$H^2(\partial_j\Xcirc_\Gamma)=\CH^*\gtimes\dR/\<u=\mm_iu\>$,
$u\in\CH^*_i$, $i=1,\ldots,n$,
and the inclusion homomorphism
$H^2(\Xcirc_\Gamma)\to H^2(\partial_j\Xcirc_\Gamma)$ is induced
by the restriction~$\res$ defined above.
\qed
\endproposition

\section{Example: pseudo-trees\label{S.trees}}

\subsection{Admissible trees and pseudo-trees}
An embedded tree $\tree\subset S^2$ is called \emph{admissible} if
all its vertices have valency~$3$ (\emph{nodes}) or~$1$
(\emph{leaves}). Each admissible tree~$\tree$ gives rise to a
skeleton~$\Gamma_\tree$: one attaches a small loop to each leaf
of~$\tree$, see Figure~\ref{fig.tree}, left.
A skeleton obtained in this
way is called a \emph{pseudo-tree}. Clearly, each
pseudo-\allowbreak tree is a
skeleton of genus~$0$.

\midinsert
\centerline{\picture{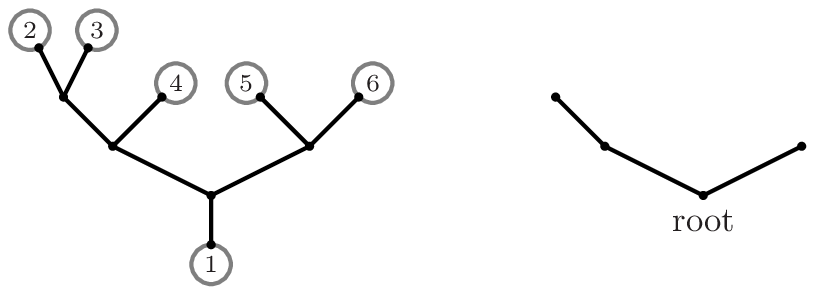}}
\figure
An admissible tree~$\tree$ (black) and skeleton $\Gamma_\tree$ (left);
the related binary tree (right)
\endfigure\label{fig.tree}
\endinsert

\paragraph\label{s.loose.end}
A nonempty admissible tree~$\tree$ has an even number $2k\ge2$ of
vertices, of which $(k-1)$ are nodes and $(k+1)$ are
leaves.
Unless $k=1$, each leaf is adjacent to a unique
node. A \emph{loose end} is a leaf sharing the same
node with an even number of other leaves. (If $k>2$, a
loose end is the only leaf adjacent to a node.)
One has
$$
\#\{\text{loose ends of~$\tree$}\}=(k+1)\bmod2.
\eqtag\label{eq.loose.ends}
$$
As a consequence, an admissible tree with $2k=0\bmod4$
vertices has a loose end.

\paragraph\label{s.marking}
A \emph{marking} of an admissible tree~$\tree$ is a choice of one of
its leaves~$v_1$. Given a marking, one can number all leaves
of~$\tree$ consecutively, starting from~$v_1$ and moving in the
clockwise direction (see Figure~\ref{fig.tree}, where the indices
of the leaves are shown inside the loops). Declaring the node
adjacent to~$v_1$ the root and removing all leaves, one obtains an
oriented rooted binary tree with $(k-1)$ vertices,
see Figure~\ref{fig.tree}, right;
conversely, an oriented rooted binary tree~$B$ gives rise to a unique
marked admissible tree: one attaches a leaf~$v_1$ at the root
of~$B$ and an extra leaf instead of each missing branch of~$B$. As
a consequence, the number of marked admissible trees with $2k$
vertices is given by the Catalan number $C(k-1)$. (Hence, the
number of unmarked admissible trees is bounded from below by
$C(k-1)/(k-1)$.)

\paragraph\label{s.distance}
The \emph{vertex distance}~$m_i$ between two consecutive
leaves~$v_i$, $v_{i+1}$ of a marked admissible
tree~$\tree$ is the vertex length of the shortest left turn
path
in~$\tree$ from~$v_i$ to~$v_{i+1}$. For example,
in Figure~\ref{fig.tree} one has
$(m_1,m_2,m_3,m_4,m_5)=(5,3,4,5,3)$. The vertex distance between
two leaves~$v_i$, $v_j$, $j>i$, is defined to be
$\sum_{k=i}^{j-1}m_k$; it is the vertex length of the shortest
left turn path connecting~$v_i$ to~$v_j$ {\em in the associated
skeleton~$\Gamma_\tree$},
\cf. Figure~\ref{fig.paths}, left, in Subsection~\ref{s.all.A} below.

\paragraph\label{s.QT}
Given a marked admissible tree~$\tree$
with $2k$ vertices, define an integral lattice $\CQ_\tree$
as follows: as a group, $\CQ_\tree$ is freely generated by $k$ vectors
$\bq_i$, $i=1,\ldots,k$ (informally corresponding to pairs
$(v_i,v_{i+1})$ of consecutive leaves),
and the products are given by
$$
\bq_i^2=m_i-2,\qquad
\bq_i\cdot\bq_j=1\ \text{if $|i-j|=1$},\qquad
\bq_i\cdot\bq_j=0\ \text{if $|i-j|\ge2$},
$$
where $m_i$, $i=1,\ldots,k$,
is the vertex distance from~$v_i$ to~$v_{i+1}$.
Next, define the \emph{characteristic functional}
$$
\chi_\tree:=\sum_{i=1}^km_i\bq_i^*\in\CQ_\tree^*.
\eqtag\label{eq.chi}
$$

\subsection{Contractions}
An \emph{elementary contraction} of an admissible tree~$\tree$ is
a new admissible tree~$\tree'$ obtained from~$\tree$ by removing two
leaves adjacent to the same node (and thus converting this node to
a leaf), see Figure~\ref{fig.contraction}.
If $\tree$ is marked, we require in addition that the two
leaves removed should be consecutive. (In other words, we do not
allow the removal of the pair $v_{k+1}$, $v_1$.) The contraction
retains a marking: if the leaves removed are~$v_1$, $v_2$, we
assign index~$1$ to their common node, becoming a leaf; otherwise,
$v_1$ remains the first leaf in~$\tree'$.

\midinsert
\centerline{\picture{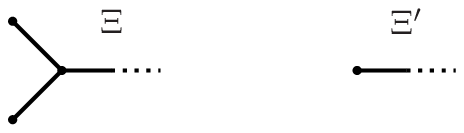}}
\figure
A tree~$\tree$ and its elementary contraction~$\tree'$\label{fig.T'}
\endfigure\label{fig.contraction}
\endinsert

By a sequence of elementary contractions any (marked) admissible
tree~$\tree$ can be reduced to a simplest tree~$\tree_0$ with two
vertices. (For proof, it suffices to consider an extremal node of
the associated binary tree: it is adjacent to two consecutive
leaves.) The resulting tree~$\tree_0$ can be identified with an
induced subtree of~$\tree$, and the reduction procedure is called a
\emph{contraction of~$\tree$ towards~$\tree_0$}. If $\tree_0$ contains a
leaf~$w$ of the original tree~$\tree$, we will also speak about a
contraction of~$\tree$ \emph{towards~$w$}.
The argument above shows that any marked admissible tree~$\tree$ can
be contracted towards its first leaf~$v_1$; similarly, $\tree$ can
be contracted towards its last leaf~$v_{k+1}$.
(In general, a contraction is \emph{not} uniquely determined by
its terminal subtree $\tree_0\subset\tree$.)

\lemma\label{contraction}
Any contraction of a marked admissible tree~$\tree$
with $2k$ vertices gives rise to an
isomorphism $\CQ_\tree\cong\CW_k$, see~\eqref{eq.UV}.
\endlemma

\proof
First,
change the sign of each even generator~$\bq_{2i}$ so that the
nontrivial exdiagonal entries of the Gram matrix
of $\CQ_\tree$ become $-1$ rather than~$1$.
The new form is represented by the graph
$$
\DD
\bul{m_1-2}\-\bul{m_2-2}\=\,\cdots\,\=\bul{m_k-2}
\eqtag\label{eq.linear.tree}
$$
where, as usual, generators are represented by the vertices (their
squares being the weights indicated) and the product of two
generator connected by an edge is~$-1$, whereas the generators not
connected are orthogonal. Whenever a graph as above has a vertex
of weight~$1$, it can be `contracted' as follows:
$$
\DD
\cdots\,\=\bul{m}\=\bul1\=\bul{n}\=\,\cdots
\quad\longmapsto\quad
\cdots\,\=\bul{m-1}\-\bul{n-1}\=\,\cdots
$$
Arithmetically, this procedure corresponds to splitting the
corresponding generator of square~$1$ as a direct summand
(passing from $\bq_{i-1}$, $\bq_i$, $\bq_{i+1}$ to
$\bq_{i-1}+\bq_i$, $\bq_i$, $\bq_{i+1}+\bq_i$,
disregarding~$\bq_i$, and leaving other generators unchanged).
On the other hand, two leaves $v_i$, $v_{i+1}$ at a vertex
distance $m_i=3$ are adjacent to the same node, and the procedure
just described establishes an isomorphism
$\CQ_\tree\cong\Z\bq_i\oplus\CQ_{\tree'}$, $\bq_i^2=1$, where $\tree'$
is the corresponding elementary contraction of~$\tree$. (In~$\tree'$,
the vertex distances just next to~$m_i$ decrease by~$1$.)
Contracting~$\tree$ to a two vertex
tree $\tree_0\subset\tree$ and
observing that $\CQ_{\tree_0}=\Z\bq_1$, $\bq_1^2=0$, one obtains an
isomorphism as in the statement.
\endproof

\Remark
Analyzing the proof, one can easily conclude that the converse of
Lemma~\ref{contraction} also holds: the lattice represented by a
linear tree~\eqref{eq.linear.tree} is isomorphic to~$\CW_k$ if and
only if,
up to the signs of the generators,
it has the form $\CQ_\tree$ for some marked admissible
tree~$\tree$.
\endRemark

According to Lemma~\ref{contraction},
a contraction
$\tree\contraction\tree_0$
sends each linear functional
$\Gf\in\CQ_\tree^*$ to a functional
$\bar\Gf\in\CW_k^*$; we will say that $\Gf$ \emph{contracts}
to~$\bar\Gf$. The following statement is straightforward.

\lemma\label{contraction.f}
If a marked admissible tree~$\tree$ with $2k$ vertices
is contracted towards
its first leaf~$v_1$, the functional
$\bq_1^*$ contracts to~$\bw^*$. If $\tree$ is contracted towards
its last
leaf~$v_{k+1}$, the functional $\bq_k^*$ contracts to
$(-1)^{k+1}\bw^*$.
\qed
\endlemma

\lemma\label{iso}
Up to isomorphism, the lattice
$\Ker\chi_\tree\subset\CQ_\tree$
does not depend on the choice of a marking of~$\tree$.
\endlemma

We postpone the proof of this lemma till next
subsection, see~\ref{pf.iso}, where a simple geometric argument is
given.

\lemma\label{contraction.2s}
If $k=2s$ is even,
the characteristic functional $\chi_\tree$, see~\eqref{eq.chi},
of a marked admissible tree~$\tree$ with $2k$ vertices
contracts\mnote{contracts} to
$$
\bar\chi=3\bv_1^*+\ldots+3\bv_s^*+\bv_{s+1}^*+\ldots+\bv_{k-1}^*
$$
\rom(up to reordering and changing the signs of the
generators~$\bv_i$\rom).
\endlemma

\proof
\Apriori, the result of contraction may depend on the choice
of\mnote{one `of' removed} a marking of~$\tree$
and on the contraction used (\cf. Remark~\ref{rem.k.odd} below).
However, we assert that, if one set of choices results in the
functional $\bar\chi$ given in the statement, then so does any
other set (up to reordering and changing the signs).
Indeed,
the divisibility of $\bar\chi$\mnote{condition added; paragraph
edited a bit} (the maximal integer $r\in\Z_{>0}$
such that $\bar\chi/r$
still takes values in~$\Z$) is the same as that of $\chi_\tree$,
and
one can easily see that, up to a scalar multiple,
$\bar\chi$ is the only functional with the
following properties:\mnote{$\bar\Gf\to\bar\chi$ everywhere}
\roster
\item\local{f1}
$\Rad\Ker\bar\chi\ne0$,
\item\local{f2}
$\det(\Ker\bar\chi/\Rad)=5k-1$, and
\item\local{f3}
the maximal root system contained in $\Ker\bar\chi/\Rad$ is
$\bA_{s-1}\oplus\bA_{s-2}$,
\endDashes
and it remains to apply Lemma~\ref{iso}. (Indeed, if
$\bar\chi=\sum_ir_i\bv_i^*+t\bw^*$ with~$t$ and all~$r_i$ coprime,
then \loccit{f1}
means that
$t=0$,
\loccit{f2} is equivalent to $\sum_ir_i^2=5k-1$, and
\loccit{f3} means that the absolute values~$|r_i|$
assume exactly
two distinct values, one $s$-fold and one $(s-1)$-fold.)

\midinsert
\centerline{\picture{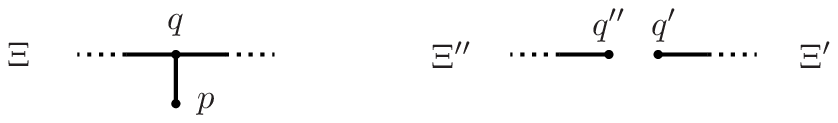}}
\figure\label{fig.loose.end}
Cutting a tree~$\tree$ at a loose end~$p$
\endfigure
\endinsert

Now, we prove the statement by induction in~$k$. For the only tree
with $4$ vertices (the case $k=2$) it is straightforward. Consider
a tree~$\tree$ with $4s\ge8$ vertices. In view
of~\eqref{eq.loose.ends}, $\tree$ has a loose end~$p$, which is the only
leaf adjacent to a certain node~$q$. Remove~$p$ and double~$q$,
cutting~$\tree$ into two trees~$\tree'$ and~$\tree''$ containing the
copies~$q'$ and~$q''$ of~$q$, respectively,
see Figure~\ref{fig.loose.end}. We may assume that
$\tree'$ contains no loose ends of the original tree~$\tree$, as otherwise
we could use that extra loose end instead of~$p$. Then,
$q'$\mnote{$q'$}
is
the only loose end of~$\tree'$ and, due to~\eqref{eq.loose.ends}, the
number of vertices in~$\tree'$ is $4s'=0\bmod4$. By additivity, the
number of vertices in~$\tree''$ is $4(s-s')=4s''=0\bmod4$. If
necessary, interchange~$\tree'$ and~$\tree''$ so that $\tree'$ is to the right
from~$p$, as in Figure~\ref{fig.loose.end}, and mark the trees
so that $q'=v'_{2s'+1}$ is the last leaf of~$\tree'$ and
$q''=v''_1$ is the first leaf of~$\tree''$. Then,
mark~$\tree$ so that $v_1=v'_1$.

Contract~$\tree'$ and~$\tree''$ towards~$q'$ and~$q''$, respectively.
This
procedure contracts~$\tree$ to a tree with a single node~$q$.
Disregarding the generators~$\bv_i'$ and $\bv_j''$
that are split off\mnote{`that are split off'}
during the contraction (in the obvious sense, they
are the same for~$\tree$ and~$\tree'$, $\tree''$), one arrives at the
quadratic form $\Z\bw'\oplus\Z\bw''$, $(\bw')^2=(\bw'')^2=1$,
$\bw'\cdot\bw''=-1$. Here, the squares of the generators
resulting from~$\tree$ differ by~$1$ from those resulting from~$\tree'$
and~$\tree''$, as so do the corresponding
vertex distances. For the same
reason, the
characteristic functional $\chi_\tree$
can be
identified with
$\chi_{\tree'}+(\bq_{s'}')^*+\chi_{\tree''}+(\bq''_1)^*$.
Due to the
induction hypothesis and Lemma~\ref{contraction.f}, it contracts
(in the obvious notation) to
$\bar\chi'-(\bw')^*+\bar\chi''+(\bw'')^*$,
and one last contraction gives the
statement for~$\tree$.
\endproof

As a corollary, we get a partial result for the case of $k$~odd.

\lemma\label{contraction.2s-1}
If $k=2s-1$ is odd and
a marked admissible tree~$\tree$ with $2k$ vertices is contracted
towards its last leaf $v_{k+1}$, the functional
$\chi_\tree$ contracts to
$$
\bar\chi=3\bv_1^*+\ldots+3\bv_{s-1}^*+\bv_s^*+\ldots+\bv_{k-1}^*+2\bw^*
$$
\rom(up to reordering and changing the signs of
the generators~$\bv_i$\rom).
\endlemma

\proof
Convert~$v_{k+1}$ to a node by attaching two extra leaves,
contract the resulting tree~$\tree'$ with $4s$ vertices towards its
last leaf, apply Lemma~\ref{contraction.2s}, and use
Lemma \ref{contraction.f} to compensate for the difference
between~$\tree$ and~$\tree'$.
\endproof

\Remark\label{rem.k.odd}
In the case of $k=2s-1$ odd, the resulting functional~$\bar\chi$
\emph{does}
depend on the choice of a contraction used.
\endRemark

\corollary\label{coprime}
If $\Gamma$ is a marked pseudo-tree with $2k\ge6$ vertices, then
the vertex distances $m_i$ are coprime\rom:
$\operatorname{g.c.d.}(m_1,\ldots,m_k)=1$.
\qed
\endcorollary

\subsection{The case of all loops of type~$\tA_0^*$}\label{s.all.A}
Consider a pseudo-tree $\Gamma=\Gamma_\tree$ and
choose the homological invariant so
that the singular fibers inside the loops attached to~$\tree$
are all of
type~$\tA_0^*$. (This choice corresponds to the boundary
orientation of each edge bounding a loop: if $v_i$ is a leaf and
$\m1(v_i)$ belongs to the original tree~$\tree$, then $\m2(v_i)$
is the tail of the new edge attached at~$v_i$.\mnote{explanation
added}
The orientations of the
edges of the original tree are irrelevant.)
Then the fiber inside the outer region
of~$\Gamma$ is of type $\tA_{5k-2}$ if $k$ is even or
$\tD_{5k+3}$ if $k$ is odd.

Pick a marking of~$\tree$, see~\ref{s.marking}, and let
$n_i=\sum_{j=i}^km_j$, $i=1,\ldots,k$, be the vertex distance
from~$v_i$ to~$v_{k+1}$, see~\ref{s.distance}.
In the computation below, we retain the notation $a,b,c\in\CH$ for
the three special elements of~$\CH$ introduced
in~\ref{s.tripods}.\mnote{a reminder added}

Mark~$\Gamma$
at each leaf~$v_i$ so that $\m1(v_i)$ belongs to the original
tree~$\tree$, see~\ref{s.orientation}. Let~$\xi_i$ be the boundary
of the loop attached at~$v_i$, and denote by $\CH^0_\Gamma$ the
subgroup spanned by the classes $\fc{\xi_i,a}$, $i=1,\ldots,k+1$.
One has $\CH^0_\Gamma\subset\Rad\CH_\Gamma$, \cf.
Example~\ref{ex.A}.
Taking into account\mnote{expanded}
constraints~\iref{rem.1}{d.vertex.alt} at~$v_i$
and~\iref{s.tripods}{d.edge} at~$\xi_i$, one
concludes that the restriction of each
element $h\in\CH_\Gamma$ to the three ends constituting~$v_i$ is a
linear combination of
$a\gtimes\m2(v_i)-b\gtimes\m3(v_i)=\fc{\xi_i,a}\in\CH^0_\Gamma$
and the element
$$
c\gtimes\m1(v_i)+b\gtimes\m2(v_i)+a\gtimes\m3(v_i).
\eqtag\label{eq.restr}
$$
Hence,
modulo~$\CH^0_\Gamma$
this restriction is a multiple of~\eqref{eq.restr},
and a dimension
count using Corollary~\ref{independent} and
Proposition~\ref{kernel}
shows that each
linear combination of elements~\eqref{eq.restr}, $i=1,\ldots,k+1$,
extends to an element of
$\CH_\Gamma/\CH^0_\Gamma$ in at most one way.
To find a simpler basis,
consider the subgroup of $\CH_\Gamma$ consisting of the vectors
satisfying all but one conditions~\iref{s.tripods}{d.edge}
and~\iref{rem.1}{d.vertex.alt}:
namely, relax~\iref{rem.1}{d.vertex.alt} at~$v_{k+1}$ to
$$
h_{\m1(v_{k+1})}+\X h_{\m2(v_{k+1})}+\X^2h_{\m3(v_{k+1})}=0\bmod b.
\eqtag\label{eq.relaxed}
$$
Let $\CH'_\Gamma$ be the quotient of this subgroup by
$\CH^0_\Gamma$. It is freely generated by the elements
$$
\be_i:=\Ge_ib\gtimes\m2(v_i)+\Ge_ia\gtimes\m3(v_i)
 +\pp{\Gg_i,\Ge_ic}
 +b\gtimes\m2(v_{k+1})+a\gtimes\m3(v_{k+1}),
$$
$i=1,\ldots,k$, where $\Gg_i$ is the shortest left turn path
from~$v_i$ to~$v_{k+1}$ and
the signs $\Ge_i=\pm1$ are chosen so that the
monodromy\mnote{monodromy spelled out for convenience}
$$
\mm_{\Gg_i}=\pm\Y(\X\Y)^{n_i-2}=
 \pm\bmatrix0&-1\\1&n_i-2\endbmatrix
$$
take $\Ge_ic$ to the
element $u_i:=c+n_ib$; these signs depend on the orientations of
the edges of the original tree~$\tree$.
Informally,\mnote{another explanation added}
$\be_i$ is obtained by extending~\eqref{eq.restr} along~$\Gg_i$,
see Definition~\ref{def.monodromy}, and `closing' it at~$v_{k+1}$
to satisfy the relaxed set of conditions;
condition~\eqref{eq.relaxed} was chosen so that the latter closure
exists and is unique modulo~$\CH^0_\Gamma$:
one merely disregards the term~$n_ib$
in $u_i$ above
and completes $c\gtimes\m1(v_{k+1})$ to~\eqref{eq.restr}.
The supports of~$\be_i$
are shown in shades of grey in Figure~\ref{fig.paths},
left; after a shift, they can be made pairwise disjoint except in a
neighborhood of the last vertex~$v_{k+1}$.

\midinsert
\centerline{\picture{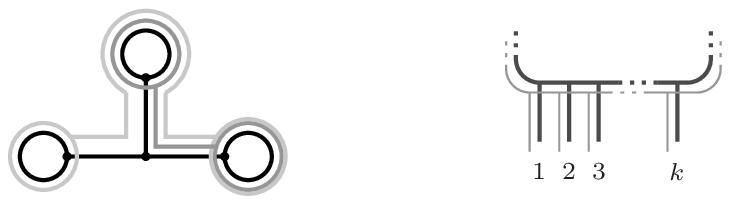}}
\figure
Supports of $\be_i$ (left) and their shifts (right)
\endfigure\label{fig.paths}
\endinsert

Bringing back\mnote{`bringing back'}
the last relation~\iref{rem.1}{d.vertex.alt}
at~$v_{k+1}$, one can see
that the subgroup $\CH_\Gamma/\CH^0_\Gamma\subset\CH'_\Gamma$
is the kernel $\Ker\Gf$, where
$$
\Gf=\sum_in_i\be_i^*.
$$
(The multiples of $n_ib_i$ disregarded in the
construction of~$\be_i$ must sum up to zero.)\mnote{more
explanation}
The self-intersection of a cycle $\sum_ir_i\be_i\in\Ker\Gf$
(assuming that it \emph{is} a cycle)
can be computed geometrically,
by shifting all paths `to the left'; it is given by
$$
\Bigl(\sum_ir_i\be_i\Bigr)^2=
 -\sum_ir_i^2-\Bigl(\sum_ir_i\Bigr)^2
 -\sum_{1\le i<j\le k}r_ir_j(u_i\cdot u_j),
\eqtag\label{eq.square}
$$
where $u_i\cdot u_j=n_i-n_j$. During the shift, the supports can
be kept pairwise disjoint except in a
small neighborhood~$U$ of~$v_{k+1}$; the
shift inside~$U$ is shown in light solid lines in
Figure~\ref{fig.paths}, right.\mnote{explanation added}
The $i$-th term of the first sum in~\eqref{eq.square}
is the contribution of the
self-intersection of $r_i\be_i$ in a neighborhood of~$v_i$,
\cf.~\eqref{eq.form.alt}. The last two terms are contributed
by~$U$. To compute this contribution, one should bring all
$1$-cycles
in the fibers
to the same basis (\cf. the proof of Theorem~\ref{tripods}); we
choose the basis in $\CH_{\m1(v_{k+1})}$. Then the $1$-cycle over the
$i$-th vertical segment in Figure~\ref{fig.paths}, right, is~$u_i$.
The $1$-cycle over the left arced segment is
$w_0:=\X^2(ra)=rb$, where $r=\sum_ir_i$, and the
$1$-cycles over the
consecutive (left to right) horizontal segments, concluding with
the right arced segment, are
$w_i:=w_0+\sum_{j=1}^ir_ju_j$, $i=1,\ldots,k$.
(Recall that $\sum_ir_i\be_i$ is \emph{assumed} a cycle.) The
intersection points are all seen in the figure, and the total
contribution from~$U$ is $-\sum_{i=1}^kw_i\cdot u_i$, which
simplifies to the last two terms in~\eqref{eq.square}.

Since we are only interested in the values
of~\eqref{eq.square} on the kernel $\Ker\Gf$,
we can add
to~\eqref{eq.square} the quadratic expression
$$
\Bigl(\sum_ir_i\Bigr)\Bigl(\sum_in_ir_i\Bigr).
$$
Now, extend the new quadratic form
to the whole group $\CH'_\Gamma$ and consider the corresponding
symmetric bilinear form; in the basis $\{\be_1,\ldots,\be_k\}$ it
is given by the matrix
$E=[e_{ij}]$, where
$e_{ii}=n_i-2$ and $e_{ij}=n_{\max\{i,j\}}-1$ for $i\ne j$.
It is straightforward that, for $i<j<k$, one has
$$
(\be_i-\be_j)\cdot\be_k=0,\quad
(\be_i-\be_j)\cdot\be_j=1,\quad
(\be_i-\be_j)\cdot\be_i=n_i-n_j.
$$
Hence, in the new basis $\bq_i=\be_i-\be_{i+1}$, $i=1,\ldots,k-1$,
$\bq_k=\be_k$ the form
turns into~$\CQ_\tree$, see~\ref{s.QT}, and the functional~$\Gf$
above turns into~$\chi_\tree$. Finally, there is an isomorphism
$$
\CH_\Gamma/\Rad=\Ker\chi_\tree/\Rad.
\eqtag\label{eq.CH=Ker}
$$

\subsubsection{Proof of Lemma~\ref{iso}}\label{pf.iso}
The statement follows from~\eqref{eq.CH=Ker} and the fact that
the left hand side does not
depend on the choice of a marking of~$\tree$.
\qed

\subsection{Proof of Theorems~\ref{th.tree.0.even}
and~\ref{th.tree.0.odd}}\label{pf.tree.0}
The skeleton~$\Gamma$
of an
extremal\mnote{`extremal'}
elliptic surface~$X$ as in the theorems is necessarily a
pseudo-tree, $\Gamma=\Gamma_\tree$,
and the singular fibers of~$X$ inside the loops
of~$\Gamma$ are all of type~$\tA_0^*$. (One has $k=2s$, $t=0$ in
Theorem~\ref{th.tree.0.even} and $k=2s-1$, $t=1$ in
Theorem~\ref{th.tree.0.odd}.)
Hence, in view of Corollary~\ref{cor.CT},
the lattice~$\CT_X$
is given by~\eqref{eq.CH=Ker}, and its structure is described by
Lemmas~\ref{contraction}, \ref{contraction.2s},
and~\ref{contraction.2s-1}. (In the case of $k=2s-1$
odd, the quotient map $\CW_k\to\CW_k/\Rad=\CV_{k-1}$ projects
$\CT_X$ to an even index~$2$
sublattice of $\CV_{k-1}$; by definition, it is
$\bD_{k-1}$.)
\qed

\subsection{The case of one type~$\tD_5$ fiber}\label{s.one.D}
Now, choose the homological invariant so that one of the loops
contain a type~$\tD_5$ fiber and mark~$\tree$ so that
this loop is attached to the last leaf~$v_{k+1}$.
Let $\xi_i$ and $\Gg_i$
be as in Subsection~\ref{s.all.A},
denote by $\CH^0_\Gamma$ the subgroup spanned by $\fc{\xi_i,a}$,
$i=1,\ldots,k$ (note that the index runs to~$k$ rather than
$k+1$), and let
$\CH''_\Gamma$ be the subgroup of
$(\CH_\Gamma/\CH^0_\Gamma)\otimes\Q$
generated over~$\Z$ by the rational cycles
$$
\gathered
\be_i:=\Ge_ib\gtimes\m2(v_i)+\Ge_ia\gtimes\m3(v_i)
 +\pp{\Gg_i,\Ge_ic}
 +v(n_i)\gtimes\m2(v_{k+1})+w(n_i)\gtimes\m3(v_{k+1}),\\
\text{where}\quad
v(n)=\frac{n}2a+\frac{2-n}4b\quad
\text{and}\quad
w(n)=\frac{2-n}4a-\frac{n}2b.
\endgathered
$$
(The vectors $v(n)$, $w(n)$ are chosen to `close' the chain
over~$v_{k+1}$, as solutions to the system
$(c+nb)+\X v+\X^2w=v+\Y w=0$.)
Then, $\CH_\Gamma/\CH^0_\Gamma$ is the
index~$4$ subgroup of~$\CH''_\Gamma$ defined by
the parity condition
$$
\psi(x)=0\bmod4,\quad
\text{where}\quad
\psi=\sum_i(n_i-2)\be_i^*.
$$
The intersection indices $\be_i\cdot\be_j$ can easily be computed
either using Theorem~\ref{tripods} or as in
Subsection~\ref{s.all.A}. One has
$$
\be_i^2=\frac14(n_i+2)(n_i-2),\qquad
\be_i\cdot\be_j=\frac14(n_i+2)(n_j-2)-1\quad
 \text{for $i<j$}.
$$
In the new basis $\bq_i=\be_i-\be_{i+1}$, $i=1,\ldots,k-1$,
$\bq_k=\be_k$
the functional~$\psi$ above takes the form
$$
\psi=\sum_im'_i\bq_i^*,\quad
\text{where $m'_i=m_i$ for $i=1,\ldots,k-1$ and $m'_k=m_k-2$},
\eqtag\label{eq.parity.m.1}
$$
and the intersection indices are
$$
\bq_i^2=\frac{(m_i')^2}4+m_i-2,\quad
\bq_i\cdot\bq_{i+1}=\frac{m_i'm_{i+1}'}4+1,\quad
\bq_i\cdot\bq_j=\frac{m_i'm_j'}4,\ j>i+1.
$$
In other words, one can identify~$\CH''_\Gamma$ with the
group~$\CQ_\tree$ supplied with the modified bilinear
form
$$
x\otimes y\mapsto x\cdot y+\frac14\psi(x)\psi(y),
$$
where $\,\cdot\,$ is the original form on~$\CQ_\tree$;
under this identification,
$\psi=\chi_\tree-2\bq_k^*$.

\subsection{Proof of Theorems~\ref{th.tree.1.even}
and~\ref{th.tree.1.odd}}\label{pf.tree.1}
The skeleton~$\Gamma$
of an
extremal\mnote{`extremal'}
elliptic surface~$X$ as in the theorems is necessarily a
pseudo-tree, $\Gamma=\Gamma_\tree$,
and the singular fibers of~$X$ inside the loops
of~$\Gamma$ are one copy of~$\tD_5$ and $k$ copies of~$\tA_0^*$.
(One has $k=2s$, $t=2$ in
Theorem~\ref{th.tree.1.even} and $k=2s-1$, $t=1$ in
Theorem~\ref{th.tree.1.odd}.)
Mark~$\tree$ as in Subsection~\ref{s.one.D} and contract it
towards its last leaf~$v_{k+1}$, establishing an isomorphism
$\CQ_\tree=\CW_k$, see Lemma~\ref{contraction}.
Due to Lemmas~\ref{contraction.2s} and~\ref{contraction.2s-1}, the
functional $\psi=\chi_\tree-2\bq_k^*$ contracts to
$$
\bar\psi=
\cases
3\bv_1^*+\ldots+3\bv_s^*+\bv_{s+1}^*+\ldots+\bv_{k-1}^*-2\bw^*,
&\text{if $k=2s$ is even},\\
3\bv_1^*+\ldots+3\bv_{s-1}^*+\bv_s^*+\ldots+\bv_{k-1}^*,
&\text{if $k=2s-1$ is odd}.
\endcases
$$
(The correction term $-2\bw^*$ is given by
Lemma~\ref{contraction.f}.)
Thus, due to
Corollary~\ref{cor.CT}
and
the results of Subsection~\ref{s.one.D},
one has
$\CT_X=\{x\in\CW'_k\,|\,\bar\psi(x)=0\bmod4\}/\Rad$,
where $\CW'_k$ is
$\CW_k$ with the modified bilinear form
$x\otimes y\mapsto x\cdot y+\frac14\bar\psi(x)\bar\psi(y)$.

If $k$ is odd, the kernel
$\Rad\CT_X$ is generated by~$\bw$, and
passing to the quotients
$\CT_X/\bw\subset\CW'_k/\bw$
one obtains the description given in
Theorem~\ref{th.tree.1.odd}.

If $k$ is even,
one has an
orthogonal decomposition $\CT_X=\Ker\bar\psi\oplus\Z\bx$, where
$\bx=2\bw$,
and $\Ker\bar\psi$ is generated by the vectors
$\bv_s-\bv_{s+1}+\bw$, $\bv_1+\bv_2+3\bw$, and
$\bv_i-\bv_{i+1}$, $i=1,\ldots,s-1,s+1,\ldots,k-2$.
It is
immediate
that $\Ker\bar\psi\cong\bD_{k-1}$.
\qed

\subsection{Proof of Theorem~\ref{th.group}}\label{pf.group}
We use Zariski--van Kampen's method~\cite{vanKampen} applied to
the ruling of~$\Sigma$.
The braid monodromy is computed using~\cite{degt.kplets}.

If there is a type~$\tD_5$ fiber, mark~$\tree$ as explained in
Subsection~\ref{s.one.D}; otherwise, mark it arbitrarily.
Mark~$\Gamma$ at~$v_{k+1}$ as shown in Figure~\ref{fig.braids}, so
that $\m1(v_{k+1})$ belongs to the
original tree~$\tree$. Take the fiber~$F$
over~$v_{k+1}$ for the reference fiber, and let
$\{\Ga_1,\Ga_2,\Ga_3\}$ be a canonical basis in~$F$ defined by the
chosen marking, see~\cite{degt.kplets}. Let~$\Gd_i$ be
the path in the base composed of the loop
of~$\Gamma$ at~$v_i$, $i=1,\ldots,k+1$,
connected to~$v_{k+1}$ by the shortest left
turn path {\em ending at $\m2(v_{k+1})$}, see
Figure~\ref{fig.braids}. According to~\cite{degt.kplets}, the
braid monodromy~$\mm_i$ along~$\Gd_i$ is given by
$$
\mm_i=\Gs_1^{n_i}\Gs_2\Gs_1^{-n_i},\ i=1,\ldots,k,\qquad
 \mm_{k+1}=\Gs_2(\Gs_1\Gs_2)^{3\epsilon},
\eqtag\label{eq.braids}
$$
where $\Gs_1$, $\Gs_2$ are the Artin generators of the braid
group~$\BG3$, parameters~$n_i$ are the vertex distances introduced
in Subsection~\ref{s.all.A}, and $\epsilon=0$ or~$1$ if the
singular fiber next to~$v_{k+1}$ is of type~$\tA_0^*$ or~$\tD_5$,
respectively.
Then one has
$$
\pi_1(\Sigma\sminus(C\cup E))=
 \bigl<\Ga_1,\Ga_2,\Ga_3\bigm|
 \text{$\mm_i=\id$, $i=1,\ldots,k+1$,
 $(\Ga_1\Ga_2\Ga_3)^{k+t}=1$}\bigr>,
$$
where $k$ and~$t$ are as introduced in~\ref{s.Sigma}. Here, each
\emph{braid relation} $\mm_i=\id$ is understood as the triple of
relations $\mm_i(\Ga_j)=\Ga_j$, $j=1,2,3$;
as a consequence, for each
$\Ga\in\<\Ga_1,\Ga_2,\Ga_3\>$ one has a relation
$\mm_i(\Ga)=\Ga$.\mnote{braid relations explained more}
The last relation
in the above presentation
is
called the \emph{relation at infinity}; in its presence,
the braid relation about the remaining
singular fiber in the outer region of~$\Gamma$ can be ignored.

\midinsert
\centerline{\picture{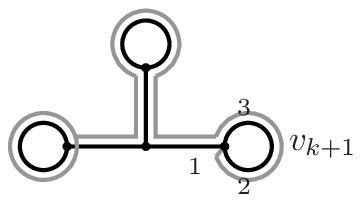}}
\figure
A loop $\Gd_i$ (grey)
\endfigure\label{fig.braids}
\endinsert

The braid relation $\mm_i(\Ga_3)=\Ga_3$, $i=1,\ldots,k$ implies
$\Ga_3=\Gs_1^{n_i}\Ga_2$.\mnote{the rest edited}
Hence one has $\Gs_1^{n_i}\Ga_2=\Gs_1^{n_j}\Ga_2$
for $1\le i,j\le k$. Since $\Gs_1$ preserves~$\Ga_3$ and the
product $\Gr:=\Ga_1\Ga_2\Ga_3$ (and $\Ga_2$, $\Ga_3$, and~$\Gr$
generate $\<\Ga_1,\Ga_2,\Ga_3\>$),
for each
$\Ga\in\<\Ga_1,\Ga_2,\Ga_3\>$
one has a relation
$\Gs_1^{n_i}\Ga=\Gs_1^{n_j}\Ga$.
Replacing~$\Ga$ with $\Gs^{-n_j}\Ga$, one can rewrite this relation
in the form $\Gs_1^{n_i-n_j}\Ga=\Ga$.

If $k>2$, the differences $n_i-n_j$, $1\le i,j\le k$, are coprime,
see Corollary~\ref{coprime}. Hence, an appropriate iteration of
the relations
$\Gs_1^{n_i-n_j}\Ga=\Ga$ obtained above results in
$\Gs_1\Ga=\Ga$, $\Ga\in\<\Ga_1,\Ga_2,\Ga_3\>$. In particular,
$\Gs_1\Ga_2=\Ga_2$, \ie., $\Ga_1=\Ga_2$. Then, the original braid
relation $\Ga_3=\Gs_1^{n_i}\Ga_2$ simplifies to $\Ga_3=\Ga_1$ or
$\Ga_3=\Ga_2$, depending on the parity of~$n_i$. In any case, one
has $\Ga_1=\Ga_2=\Ga_3$, and the group is cyclic.
\qed


\Remark\label{rem.group}
In the exceptional cases $k=1,2$, the
fundamental groups are also easily
computed. We skip details and merely indicate the result:
$$
\alignedat2
&k=2,\ t=0:\quad&&\BG3/(\Gs_1\Gs_2)^3,\\
&k=2,\ t=2:\quad&&\CG3\rtimes\CG{12},\\
&k=1,\ \text{\smash{$\tD_8$} type fiber}:\quad&&
 \Z\times\CG2,\\
&k=1,\ \text{\smash{$\tD_5$} type fiber}:\quad&&
 \Z[t]/(t^2-1)\rtimes\CG2;
\endalignedat
$$
in the last case, the generator of~$\CG2$ act on the kernel \via\
multiplication by~$t$. It follows that, for $k=1$, the trigonal
curve~$C$ is reducible.
\endRemark

\section{Generalizations\label{S.generalizations}}

In this section, we
outline
two generalizations of
Theorem~\ref{tripods}, one to surfaces with
type~$\tE$ singular fibers, and one to non-extremal surfaces.

\subsection{Extremal surfaces with $\tE$ type fibers}\label{s.E}
Let~$X$ be an extremal elliptic surface \emph{with} type~$\tE$
singular fibers. (Accidentally, at this point we can also admit
singular fibers of types $\tA_0^{**}$,
$\tA_1^*$, or~$\tA_2^*$, provided that $X$ satisfies
conditions~\iref{def.extremal}{ext.1} and~\ditto{ext.2} and has no
fibers of type~$\tD_4$.) Let $\Gamma=\Gamma_X$ be the skeleton
of~$X$; it may have \black-vertices of valency~$\le2$ or
\white-vertices of valency~$1$.
Replace these irregular
vertices with the boundaries of small disks, see
Figure~\ref{fig.e-fiber}, bottom row, converting~$\Gamma$
to a regular $3$-graph~$\Gamma'$. Unlike~$\Gamma$, the new
graph~$\Gamma'$ is a skeleton in the sense of
Definition~\ref{def.skeleton}.

\midinsert
\centerline{\picture{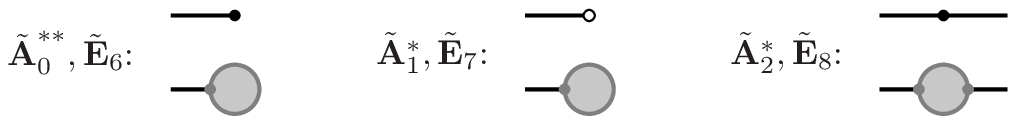}}
\figure
The modification $\Gamma'$ of $\Gamma$
\endfigure\label{fig.e-fiber}
\endinsert

Orient the new edges of~$\Gamma'$ as the boundary of the shaded
regions in Figure~\ref{fig.e-fiber}. Assign label
$\Y\in\PSL(2,\Z)$ to each edge of the original skeleton~$\Gamma$,
and label the new edges (grey in the figure) as follows:
\Dashes
\dash
type~$\tA_0^{**}$ ($\tE_6$): the label is~$\X$
(respectively,~$-\X$);
\dash
type~$\tA_1^*$ ($\tE_7$): the label is~$-\X\Y\X$
(respectively,~$\X\Y\X$);
\dash
type~$\tA_2^*$ ($\tE_8$): the two labels are either both~$\X$ or
both~$-\X$
(respectively, one label is~$\X$ and one is~$-\X$).
\endDashes
(In the last case, when two new edges are inserted,
there are two choices of the labelling;
they result in distinct homological
invariants of~$X$, \cf.~\ref{s.monodromy.E} below.)

Define $\CH\gtimes\Gamma$ to be $\CH\gtimes\Gamma'$, see
Subsection~\ref{s.tripods}, and let $\CH_\Gamma$ be the subgroup
of $\CH\gtimes\Gamma$ subject to
the following relations:
\roster
\item\local{dd.vertex}
$h_{\Ga}+\X h_{\nx\Ga}+\X^2h_{\nx^2\Ga}=0$ for each element
$\Ga\in\CE$;
\item\local{dd.edge}
$h_{e^+}+L h_{e^-}=0$ for each edge $e\in\CE\!/\!\op$
labelled~$L$.
\endroster
Similarly, let
$\CH^*\gtimes\Gamma=\CH^*\gtimes\Gamma'$
and define $\CH^*_\Gamma$ as the
quotient of $\CH^*\gtimes\Gamma$ by the subgroup spanned by the
vectors of the form
\roster
\item[3]\local{dd*.vertex}
$u\gtimes\Ga+\X^*u\gtimes(\nx\Ga)+(\X^*)^2u\gtimes(\nx^2\Ga)$
for each $u\in\CH^*$ and $\Ga\in\CE$;
\item\local{dd*.edge}
$u\gtimes e^++L^*u\gtimes e^-$ for each
$u\in\CH^*$ and each edge $e\in\CE\!/\!\op$ labelled~$L$.
\endroster
There is a natural paring
$\CH_\Gamma\otimes\CH^*_\Gamma\to\Z$.

\theorem\label{tripods.E}
There are isomorphisms $H_2(\Xcirc)=\CH_\Gamma$ and
$H^2(\Xcirc)=\CH^*_\Gamma$.
The former takes the intersection index form to the form given
by~\eqref{eq.form}\rom; the latter takes the Kronecker product to
the pairing $\CH_\Gamma\otimes\CH^*_\Gamma$ defined above.
\endtheorem

\proof
The proof repeats literally that of Theorem~\ref{tripods}: the
space $\Xcirc$ has a strict deformation retract~$X'$ which fibers
over the new graph~$\Gamma'$.
\endproof

\corollary
One has $\CT_X=\CH_\Gamma/\Rad$ and
$\Tors\MW(X)=\Tors\CH^*_\Gamma$.
\qed
\endcorollary

\paragraph\label{s.monodromy.E}
One can also mimic Definition~\ref{def.monodromy}
and define the \emph{monodromy}~$\mm_\Gg$ and
the \emph{parallel transport}
$\pp{\Gg,h_0}\in\CH\gtimes\Gamma$ along a path~$\Gg$
{\em in the new graph~$\Gamma'$}.
Part~\iref{def.monodromy}{m.edge} of the definition should be
replaced with $\mm_i=-L^{\pm1}$ for an edge $e=[\Ga_{i-1},\Ga_i]$
labelled~$L$ and $\Ga_i=e^\pm$. Under this definition, the
monodromy along the boundary of each region~$R$ other than the
shaded disks in Figure~\ref{fig.e-fiber} is still of the form
$\pm(\X\Y)^n$, where $n$ is the number of corners
(\black-vertices) in the boundary $\dR$ {\em in the original
graph~$\Gamma$}. The monodromy along the boundary of a
shaded disk~$R$ is of the form $\pm\X^{\pm1}$ or $\pm\X\Y\X\1$,
depending on the type of the singular fiber inside~$R$.
It follows that the monodromies~$\mm_{\dR}$ determine the type
specification of~$X$ and that $\mm_{\dR}$ has an invariant vector
if and only if the singular fiber inside~$R$ is stable.
In particular, one still has analogues of
Propositions~\ref{kernel} and~\ref{cohomology}.

\subsection{Non-extremal surfaces}\label{s.non-extremal}
Now, consider a Jacobian elliptic surface~$X$, not necessarily
extremal, satisfying the following conditions (\cf.
Definition~\ref{def.extremal}):
\roster
\item\local{nonext.1}
$j_X$ has no critical values other than~$0$, $1$, and~$\infty$;
\item\local{nonext.2}
each point in $j_X\1(0)$ has ramification index $(0\bmod3)$,
and each point in $j_X\1(1)$ has ramification index~$2$;
\item\local{nonext.3}
$X$ has no singular fibers of type~$\tD_4$.
\endroster
(We do not discuss whether any elliptic
surface can be deformed to one
satisfying \loccit{nonext.1}--\loccit{nonext.3}.
For each particular surface~$X$,
this can
be decided in terms of equisingular degenerations
of the dessin of~$X$, see~\cite{degt.kplets}.)

As in
Subsection~\ref{s.GammaX}, define the \emph{skeleton}
$\Gamma_X=j_X\1[0,1]$; it is a ribbon graph with all vertices of
valency $(0\bmod3)$. (The idea
of considering
skeletons with multiple
vertices rather than dessins in the sense of~\cite{degt.kplets}
was suggested to me by I.~Shimada.)
To accommodate $\Gamma=\Gamma_X$,
modify Definition~\ref{def.skeleton} by
replacing the condition $\nx^3=\id$ with the requirement that
each orbit of~$\nx$ should have length divisible by~$3$.
Then, as in Subsection~\ref{s.tripods}, define
$\CH\gtimes\Gamma=\bigoplus\CH_\Ga$, $\Ga\in\CE$, and let
$\CH_\Gamma\subset\CH\gtimes\Gamma$ be the subgroup subject to the
following conditions:
\roster
\item\local{ne.vertex}
$\sum_{i=0}^{n-1}\X^ih_{\Ga_i}=0$
for each vertex $(\Ga_0,\ldots,\Ga_{n-1})\in\CE\!/\!\nx$ of
valency~$n$;
\item\local{ne.edge}
$h_{e^+}+\Y h_{e^-}=0$ for each edge $e\in\CE\!/\!\op$.
\endroster
Also, define $\CH^*_\Gamma$ as the quotient of
$\CH^*\gtimes\Gamma$ by the image of the maps adjoint
to the left hand sides of~\loccit{ne.vertex},
\loccit{ne.edge}. There is a pairing
$\CH_\Gamma\otimes\CH^*_\Gamma\to\Z$.

Convert $\CH\gtimes\Gamma$ to a rational lattice, defining the
square~$h^2$ of
$h\in\CH\gtimes\Gamma$ to be
$$
h^2=-\sum_{\Ga\in\CE}\sum_{d=1}^{n(\Ga)-2}
 \frac{n(\Ga)-d-1}{n(\Ga)}\,h_\Ga\cdot\X^dh_{\nx^d\Ga},
\eqtag\label{eq.form.ne}
$$
where $n(\Ga)$ is the valency of the vertex represented by~$\Ga$,
\ie., the length of the orbit of~$\nx$ containing~$\Ga$. (An
alternative expression for the restriction of this form
to~$\CH_\Gamma$ is given by~\eqref{eq.form.ne.alt} below.)

\theorem\label{tripods.ne}
There are isomorphisms $H_2(\Xcirc_\Gamma)=\CH_\Gamma$ and
$H^2(\Xcirc_\Gamma)=\CH^*_\Gamma$.
The former takes the intersection index form to the form given
by~\eqref{eq.form.ne}\rom; the latter takes the Kronecker product to
the pairing $\CH_\Gamma\otimes\CH^*_\Gamma$ defined above.
\endtheorem

\proof
Again, the proof repeats literally that of Theorem~\ref{tripods}.
To compute the contribution to the self intersection~$h^2$
of a cycle
$h\in\CH_\Gamma$ by a marked $n$-valent vertex
$(\Ga_0,\ldots,\Ga_{n-1})\in\CE\!/\!\nx$, `spread out' and shift the
vertex as shown in Figure~\ref{fig.spread}. The resulting
expression is
$$
-\sum_{i=1}^{n-2}\sum_{j=0}^{i-1}\X^jh_j\cdot\X^ih_i=
 -\sum_{i=1}^{n-2}\sum_{j=0}^{i-1}h_j\cdot\X^{i-j}h_i,
\eqtag\label{eq.form.ne.alt}
$$
\cf. the proof of~\eqref{eq.square}.\mnote{ref added}
(We abbreviate $h_i=h_{\Ga_i}$ and use the fact that $\X$ is an
isometry.)
Averaging over all $n$ markings of the vertex,
one arrives at~\eqref{eq.form.ne}.
\endproof

\midinsert
\centerline{\picture{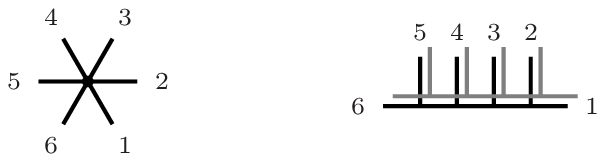}}
\figure
Spreading out and shifting a vertex
\endfigure\label{fig.spread}
\endinsert

\corollary
One has $\CT_X=\CH_\Gamma/\Rad$ and
$\Tors\MW(X)=\Tors\CH^*_\Gamma$.
\qed
\endcorollary

With the obvious modifications, the material of
Subsection~\ref{s.monodromy} extends to the general case.
One can also combine the constructions of this and the previous
subsections and consider non-extremal surfaces with $\tE$ type
singular fibers. (For the sake of simplicity, it is better to
consider a skeleton~$\Gamma$ with all \white-vertices of
valency~$\le2$ and all \black-vertices of valency either
$(0\bmod3)$ or~$\le2$.) We leave details to the reader.

\refstyle{C}

\enddocument